\documentclass{article}

\usepackage{amsfonts, amsmath, amssymb, amsthm, amscd}
\usepackage{ascmac}
\usepackage{verbatim}
\usepackage{graphics}
\usepackage[all, 2cell]{xypic}

\theoremstyle{plain}  
\newtheorem{thm}{Theorem}[section]

\theoremstyle{definition}

\newtheorem{para}[thm]{}

\theoremstyle{remark}

\DeclareMathOperator{\cA}{\mathcal{A}}
\DeclareMathOperator{\cB}{\mathcal{B}}
\DeclareMathOperator{\cC}{\mathcal{C}}
\DeclareMathOperator{\calD}{\mathcal{D}}
\DeclareMathOperator{\cE}{\mathcal{E}}
\DeclareMathOperator{\cF}{\mathcal{F}}

\DeclareMathOperator{\cI}{\mathcal{I}}

\DeclareMathOperator{\calL}{\mathcal{L}}
\DeclareMathOperator{\cM}{\mathcal{M}}
\DeclareMathOperator{\cN}{\mathcal{N}}
\DeclareMathOperator{\cO}{\mathcal{O}}
\DeclareMathOperator{\cP}{\mathcal{P}}

\DeclareMathOperator{\cS}{\mathcal{S}}
\DeclareMathOperator{\cT}{\mathcal{T}}
\DeclareMathOperator{\cU}{\mathcal{U}}

\DeclareMathOperator{\cX}{\mathcal{X}}
\DeclareMathOperator{\cY}{\mathcal{Y}}
\DeclareMathOperator{\cZ}{\mathcal{Z}}

\DeclareMathOperator{\bE}{\mathbf{E}}

\DeclareMathOperator{\bbA}{\mathbb{A}}

\DeclareMathOperator{\bbF}{\mathbb{F}}

\DeclareMathOperator{\bbK}{\mathbb{K}}

\DeclareMathOperator{\bbN}{\mathbb{N}}

\DeclareMathOperator{\bbP}{\mathbb{P}}

\DeclareMathOperator{\bbR}{\mathbb{R}}

\DeclareMathOperator{\bbZ}{\mathbb{Z}}

\DeclareMathOperator{\fF}{\mathfrak{F}}
\DeclareMathOperator{\ff}{\mathfrak{f}}

\DeclareMathOperator{\fg}{\mathfrak{g}}

\DeclareMathOperator{\fraki}{\mathfrak{i}}

\DeclareMathOperator{\fn}{\mathfrak{n}}

\DeclareMathOperator{\ft}{\mathfrak{t}}

\DeclareMathOperator{\fx}{\mathfrak{x}}

\DeclareMathOperator{\fy}{\mathfrak{y}}

\UseAllTwocells

\def\sn{\smallskip\noindent}
\def\mn{\medskip\noindent}

\def\enumidef{\renewcommand{\labelenumi}{$\mathrm{(\arabic{enumi})}$}}
\def\enumiidef{\renewcommand{\labelenumii}{$\mathrm{(\roman{enumii})}$}}

\newcommand{\add}{\operatorname{add}}

\newcommand{\Bl}{\operatorname{Bl}}

\newcommand{\cf}{\textrm{cf.}\;}
\newcommand{\Ch}{\operatorname{\bf Ch}}
\newcommand{\codim}{\operatorname{codim}}
\newcommand{\Coh}{\operatorname{\bf Coh}}
\newcommand{\coker}{\operatorname{Coker}}
\newcommand{\colim}{\operatorname{colim}}
\newcommand{\Cone}{\operatorname{Cone}}
\newcommand{\Cub}{\operatorname{\bf Cub}}

\newcommand{\dg}{\operatorname{dg}}
\newcommand{\dgCat}{\operatorname{\bf dgCat}}

\newcommand{\ExCat}{\operatorname{\bf ExCat}}

\newcommand{\gr}{\operatorname{Gr}}
\newcommand{\graded}{\operatorname{gr}}
\newcommand{\grCoh}{\operatorname{\bf GrCoh}}
\newcommand{\grMod}{\operatorname{\bf GrMod}}

\newcommand{\Ho}{\operatorname{Ho}}
\newcommand{\hocolim}{\operatorname{hocolim}}
\newcommand{\Hom}{\operatorname{Hom}}
\newcommand{\HOM}{\mathcal{H}\!\!om}
\newcommand{\Homo}{\operatorname{H}}

\newcommand{\id}{\operatorname{id}}
\newcommand{\im}{\operatorname{Im}}
\newcommand{\isoto}{\overset{\scriptstyle{\sim}}{\to}}

\newcommand{\Kos}{\operatorname{\bf Kos}}
\newcommand{\kos}{\operatorname{Kos}}

\newcommand{\Lex}{\operatorname{\bf Lex}}
\newcommand{\linf}{\leftarrowtail}
\newcommand{\loc}{\operatorname{loc}}

\newcommand{\Mod}{\operatorname{\bf Mod}}
\newcommand{\Mot}{\mathcal{M}ot}

\newcommand{\Nil}{\operatorname{Nil}}
\newcommand{\nilp}{\operatorname{nilp}}

\newcommand{\Ob}{\operatorname{Ob}}
\newcommand{\onto}[1]{\stackrel{#1}{\to}}
\newcommand{\op}{\operatorname{op}}

\newcommand{\Perf}{\operatorname{\bf Perf}}
\newcommand{\proj}{\operatorname{Proj}}
\newcommand{\projdim}{\operatorname{Projdim}}

\newcommand{\Qcoh}{\operatorname{\bf Qcoh}}
\newcommand{\qis}{\operatorname{qis}}

\newcommand{\rdef}{\twoheadrightarrow}
\newcommand{\red}{\operatorname{red}}
\newcommand{\reg}{\operatorname{reg}}
\newcommand{\rinc}{\hookrightarrow}
\newcommand{\rinf}{\rightarrowtail}

\newcommand{\Spec}{\operatorname{Spec}}
\newcommand{\simp}{\operatorname{simp}}
\newcommand{\ssm}{\smallsetminus}
\newcommand{\stab}{\operatorname{Stab}}
\newcommand{\supp}{\operatorname{Supp}}

\newcommand{\Tot}{\operatorname{Tot}}
\newcommand{\topo}{\operatorname{top}}
\newcommand{\tq}{\operatorname{tq}}
\newcommand{\typ}{\operatorname{typ}}
\newcommand{\Typ}{\operatorname{Typ}}

\newcommand{\uloc}{\operatorname{uloc}}
\newcommand{\unilp}{\operatorname{unilp}}

\title{Nilpotent invariant motives I}
\date{}

\author{Satoshi Mochizuki}

\begin{document}

\maketitle

\abstract{
The purpose of this article is to clarify the question 
what makes motives $\mathbb{A}^1$-homotopy invariance. 
}

\section*{Introduction}

The purpose of this article is to clarify the question 
what makes motives $\mathbb{A}^1$-homotopy invariance. 
We give a guide for the structure of this article. 
In section~\ref{sec:nilp inv motives}, 
we recall the construction of the stable model category of 
nilpotent invariant motives $\Mot_{\dg}^{\nilp}$ 
from \cite{Moc16b} 
and define the nilpotent invriant motives 
associated with schemes and relative exact categories. 
For a noetherian scheme $X$, there are two kind of motives associated with 
$X$ in the homotopy category $\Ho(\Mot^{\nilp}_{\dg})$, 
namely $M_{\nilp}(X)$ and $M_{\nilp}'(X)$. 
In general $M_{\nilp}(X)$ is not isomorphic to $M'_{\nilp}(X_{\red})$. 
But there exists a canonical 
isomorphism $M_{\nilp}'(X)\isoto M_{\nilp}'(X_{\red})$ and 
if $X$ is regular noetherian separated, $M(X)$ is canonically 
isomorphic to $M'(X)$. 
In section~\ref{sec:absolute geometric presentation theorem}, 
we will show absolute geometric presentation 
theorem~\ref{thm:abs geom presentation theorem}. 
Roughly saying this theorem says that 
for a regular local ring $A$ and for each non-negative integer 
$p\leq \dim A$, the topological weight $p$ part of $\Spec A$ is 
isomorpic to the Adams weight $p$ part of $\Spec A$ as nilpotent 
invariant motives. 
In section~\ref{sec:what makes motives A1homotopy invariance}, 
we establish a concept of algebraic varieties over 
(locally) noetherian abelian categories. 
By utilizing this concept we clarify the question 
what makes motives $\mathbb{A}^1$-homotopy invariance. 
We will show $\mathbb{A}^1$-homotopy invariance of $M'$ 
which does not hold for $M$ in general. 

\paragraph{Acknowledgement}
This article is the extended notes of 
my  talk given at the 
workshop ``Algebraic number theory and related topics'' 
at RIMS, Kyoto in December 2016 and 
Shinshu Topology seminar at Shinshu University 
in January 2017. 
I wish to express my deep gratitude to all 
organizers of the workshop and the seminar, 
professors Yasuo Ohno, Hiroshi Tsunogai, 
Dai Tamaki, Keiichi Sakai and Katsuhiko Kuribayashi 
for giving me the opportunity to present my work and 
I greatly appreciate Mariko Ohara for support 
during my stay in Shinshu and Tetsuya Uematsu, 
Akiyoshi Sannai and Seidai Yasuda
for stimulative 
discussion in the early stage of 
the works in 
section~\ref{sec:absolute geometric presentation theorem} and 
section~\ref{sec:what makes motives A1homotopy invariance}.

\section{Stable model category of nilpotent invariant motives}
\label{sec:nilp inv motives}

In this section, 
we review the costrcution of 
the stable model category of nilpotent invariant motives 
$\Mot_{\dg}^{\nilp}$ 
from \cite{Moc16b} and define motives associated with coherent schemes 
and relative exact categories and dash motives associated with 
noetherian schemes. 

We denote the category of small dg-categories over $\bbZ$ 
the ring of integers 
by $\dgCat$ and let 
$\Mot^{\add}_{\dg}$ and 
$\Mot^{\loc}_{\dg}$ 
be symmetric monoidal stable 
model categories of additve noncommutative motives and 
localizing motives over $\bbZ$ respectively. 
(See \cite[\S 7]{CT12}.) 
In general we denote the homotopy category of a model category $\cM$ by 
$\Ho(\cM)$. 
There are functors from $\ExCat$ the 
category of small exact categories to $\dgCat$ 
which sending a small exact category $\cE$ to 
its bounded dg-derived category $\calD_{\dg}^b(\cE)$ 
(see \cite[\S 4.4]{Kel06}.) 
and the universal functors 
$\cU_{\add}\colon\dgCat\to 
\Ho(\Mot^{\add}_{\dg})$ 
and 
$\cU_{\loc}\colon\dgCat\to 
\Ho
(\Mot^{\loc}_{\dg})$. 
We denote the compositions of these functors 
$\ExCat\to
\Ho(\Mot^{\add}_{\dg})$ and 
$\ExCat\to\Ho(\Mot^{\loc}_{\dg})$ 
by $M_{\add}$ and $M_{\loc}$ respectively. 

We start by recalling the construction of the stable model category 
of localizing non-commutative motives 
${\Mot}^{\loc}_{\dg}$ 
from \cite{Tab08}. 
First notice that 
the category $\dgCat$ carries a cofibrantly generated 
model structure whose weak equivalences are the derived Morita equivalences. 
\cite[Th\'eor\`eme 5.3]{Tab05}. 
We fix on $\dgCat$ a fibrant resolution functor $R$, 
a cofibrant resolution functor $Q$ and a left framing 
$\Gamma_{\ast}$ (see \cite[Definition 5.2.7 and Theorem 5.2.8]{Hov99}.) 
and we also fix a small full subcategory 
$\dgCat_f\hookrightarrow\dgCat$ 
such that it contains all finite dg cells and 
any objects in $\dgCat_f$ are ($\bbZ$-)flat 
and homotopically finitely presented (see \cite[Definition 2.1 (3)]{TV07}.) 
and $\dgCat_f$ is closed uder the operations $Q$, $QR$ and 
$\Gamma_{\ast}$ and $\otimes$. 
The construction below does not depend upon a choice of $\dgCat_f$ 
up to Dwyer-Kan equivalences. 
Let $s\widehat{\dgCat_f}$ and 
$s\widehat{{\dgCat_f}_{\bullet}}$ be the category of simplicial presheaves and that of pointed simplicial 
preshaves on $\dgCat_f$ respectively. 
We have the projective model structures on $s\widehat{\dgCat_f}$ 
and $s\widehat{{\dgCat_f}_{\bullet}}$ 
where the weak equivalences and the fibrations are 
the termwise simplicial weak equivalences and termwise Kan fibrations 
respectively. (see \cite[Theorem 11.6.1]{Hir03}.) 
We denote the class of derived Morita equivalences in 
$\dgCat_f$ by $\Sigma$ and 
we also write $\Sigma_{+}$ for the image of $\Sigma$ 
by the composition of the Yoneda embedding 
$h\colon \dgCat_f \to s\widehat{\dgCat}$ and 
the canonical functor 
${(-)}_{+}\colon s\widehat{\dgCat_f}\to 
s\widehat{{\dgCat_f}_{\bullet}}$. 
Let $P$ be the canonical map $\emptyset \to h(\emptyset)$ in 
$s\widehat{\dgCat_f}$ and we write $P_{+}$ 
for the image of $P$ by the functor ${(-)}_{+}$. 
We write $L_{\Sigma,P}s\widehat{{\dgCat_f}_{\bullet}}$ 
for the left Bosufield localization of 
$s\widehat{{\dgCat_f}_{\bullet}}$ by the 
set $\Sigma_{+}\cup\{P_{+}\}$. 
The Yoneda embedding functor induces a functor 
$$\bbR\underline{h}\colon \Ho(\dgCat)\to\Ho(L_{\Sigma,P}s\widehat{{\dgCat_f}_{\bullet}})$$
which associates any dg category $\cA$ to 
the pointed simplicial presheaves on $\dgCat_f$:
$$\bbR\underline{h}(\cA)\colon\cB\mapsto 
{\Hom(\Gamma_{\ast}(Q\cB),R(\cA))}_{+}.$$
Let $\cE$ be the class of morphisms in 
$L_{\Sigma,P}s\widehat{{\dgCat_f}_{\bullet}}$ of shape 
$$\operatorname{Cone}[\bbR\underline{h}(\cA)\to 
\bbR\underline{h}(\cB)]\to \bbR\underline{h}(\cC)$$
associated to each exact sequence of dg categories 
$\cA\to\cB\to\cC$, 
with $\cB$ in $\dgCat_f$ 
where $\Cone$ means homotopy cofiber. 
We write ${\Mot}_{\dg}^{\uloc}$ 
for the left Bosufield localization of 
$L_{\Sigma,P}s\widehat{{\dgCat_f}_{\bullet}}$ by $\cE$ 
and call it the {\it model category of 
unstable localizing non-commutative motives}. 
Finally we write ${\Mot}_{\dg}^{\loc}$ 
for the stable symmetric monoidal model 
category of symmetric $S^1\otimes 1\!\!\!1$-spectra on 
${\Mot}_{\dg}^{\uloc}$ 
(see \cite[\S 7]{Hov01}.) and call it 
the {\it model category of localizing non-commutative motives}. 
Next we construct the stable model category ${\Mot}_{\dg}^{\nilp}$. 
First recall that 
we say that a non-empty full subcategory $\cY$ of 
a Quillen exact category $\cX$ 
is a {\it topologizing subcategory} of $\cX$ if 
$\cY$ is closed under finite direct sums and closed under 
admissible sub- and quotient objects. 
The naming of the term `topologizing' comes from 
noncommutative geometry of abelian categories 
by Rosenberg. (See \cite[Lecture 2 1.1]{Ros08}.) 
We say that a full subcategory $\cY$ of 
an exact category $\cX$ is a {\it Serre subcategory} 
if it is an extensional closed topologizing subcategory of $\cX$. 
For any full subcategory $\cZ$ of $\cX$, 
we write ${}^S\!\!\!\sqrt{\cZ}$ 
for intersection of all Serre subcategories 
which contain $\cZ$ and call it the 
{\it Serre radical of $\cZ$} 
({\it in $\cX$}). 
We say that an object $x$ in an exact category is {\it noetherian} 
if any ascending filtration of admissible subobjects of $x$ is stational. 
We say that an exact category $\cE$ is {\it noetherian} if 
it is essentially small and 
all objects in $\cE$ 
are noetherian.

\begin{para}
\label{df:nilp immersion}
{\bf Definition (Nilpotent immersion).}\ \ 
Let $\cA$ 
be a noetherian abelian category and let 
$\cB$ a topologizing subcategory. 
We say that {\it $\cB$ satisfies the d\'evissage condition} 
({\it in $\cA$}) or say that the inclusion 
$\cB\rinc \cA$ 
is a {\it nilpotent immersion}
if one of the following equivalent conditions holds:
\begin{enumerate}
\item
For any object $x$ in $\cA$, 
there exists a finite filtration of 
monomorphisms 
$$x=x_0\leftarrowtail x_1\leftarrowtail x_2\leftarrowtail \cdots 
\leftarrowtail x_n=0$$
such that for every $i<n$, $x_i/x_{i+1}$ is 
isomorphic to an object in $\cB$. 

\item 
We have the equality 
$$\cA={}^S\!\!\!\sqrt{\cB}.$$
\end{enumerate}
(For the proof of the equivalence of the conditions above, 
see \cite[3.1]{Her97}, \cite[2.2]{Gar09}.) 
\end{para}

\begin{para}
\label{df:nilp motives}
{\bf Definition (Nilpotent invariant motives).}\ \ 
We write $\cN$ for the class of morphisms in 
${\Mot}_{\dg}^{\uloc}$ of shape 
$$\bbR\underline{h}(\cC\otimes\calD_{\dg}^b(\cB))\to\bbR\underline{h}
(\cC\otimes\calD_{\dg}^b(\cA))$$
associated with each noetherian abelian category $\cA$ 
and each nilpotent immersion $\cB\rinc \cA$ and each 
small dg-category $\cC$ in $\dgCat_f$. 
We write ${\Mot}_{\dg}^{\unilp}$ 
the left Bousfield localization of 
${\Mot}_{\dg}^{\uloc}$ by 
$\cN$ and 
call it the 
{\it model category of unstable nilpotent invariant non-commutative motives}. 
This category naturally becomes a symmetric monoidal model category 
by \cite[Theorem 5.7.]{CT12}
Finally we write ${\Mot}_{\dg}^{\nilp}$ 
for the stable model 
category of symmetric $S^1\otimes 1\!\!\!1$-spectra on 
${\Mot}_{\dg}^{\unilp}$ 
and call it 
the {\it stable model category of nilpotent invariant non-commutative motives}. 
We denote the compositions of the following functors
$$
\dgCat\to\Ho(\dgCat)
\onto{\bbR\underline{h}} 
\Ho(L_{\Sigma,P}s\widehat{{{\dgCat}_f}_{\bullet}})\to
\Ho({\Mot}_{\dg}^{\unilp})
\onto{\Sigma^{\infty}}
\Ho({\Mot}_{\dg}^{\nilp})
$$
by $\cU_{\nilp}$ and we write $M_{\nilp}$ for the compositions of 
the following functors
$$\ExCat\onto{\calD^b_{\dg}}
\dgCat\onto{\cU_{\nilp}}\Ho({\Mot}_{\dg}^{\nilp}).$$
\end{para}

We recall the conventions of relative exact categories 
from \cite{Moc13b}. 
A {\it relative exact category} $\bE=(\cE,w)$ is a pair of 
exact category $\cE$ and a class of morphisms $w$ in $\cE$. 
We call $\cE$ and $w$ the {\it underlying exact category of $\bE$} 
and the {\it class of weak equivalences of $\bE$} and 
denote it by $\cE_{\bE}$ and $w_{\bE}$ respectively. 
We say that a relative exact category is {\it small} if 
its underlying exact category is small. 
A {\it relative exact functor} $f\colon \bE \to \bE'$ 
is an exact functor $f\colon\cE_{\bE} \to \cE_{\bE'}$ such that 
$f(w_{\bE})\subset w_{\bE'}$. 
We say that a relative exact category is {\it extensional} 
if its class of weak equivalences satisfies the extension axiom in 
\cite{Wal85}. 
We say that a relative exact category $\bE$ is {\it biWaldhausen} 
if both $w_{\bE}$ and $w^{\op}_{\bE}$ satisfy the gluing 
axiom in $\cE_{\bE}$ and $\cE_{\bE}^{\op}$ respectively. 
For a relative exact category $(\cE,w)$, 
let $\cE^w$ be the full subcategory of $\cE$ consisting of 
those object $x$ such that the canonical morphism $0\to x$ is in $w$. 
If we assume $(\cE,w)$ is either extensional or biWaldhausen, then 
$\cE^w$ is an exact category such that 
the exact functor 
$\cE^w\rinc \cE$ is exact and reflects exactness. 
(See \cite[Proposition 2.4.]{Moc13b}.)

\begin{para}
\label{df:motives ass with rel exact cat}
{\bf Definition (Motives associated with relative exact categories).}\ \ 
Let $\bE=(\cE,w)$ be a small 
relative exact category which is either extensional or 
biWaldhausen. Then 
we define $\calD^b_{\dg}(\bE)=\calD^b_{\dg}(\cE,w)$ to be a dg-category 
by setting $\calD^b_{\dg}(\bE):=\calD^b_{\dg}(\cE)/\calD^b_{\dg}(\cE^w)$ 
the Drinfeld quotient of $\calD^b_{\dg}(\cE)$ by $\calD^b_{\dg}(\cE^w)$ and 
we call it the {\it bounded dg-derived category of $\bE$}. 
For $\#\in\{\add,\loc,\nilp\}$, we set 
$M_{\#}(\bE)=M_{\#}(\cE,w):=\cU_{\#}(\calD^b_{\dg}(\bE))$ 
and call it {\it additive} (resp. {\it localizing, nilpotent invariant}) 
{\it motives associated with $\bE$}. 
For a small exact category $\cE$, we write $i_{\cE}$ for the class of 
all isomorphisms in $\cE$. 
Then we have the canonical identification 
$M_{\#}(\cE,i_{\cE})=M_{\#}(\cE)$. 
\end{para}

A scheme $X$ is {\it coherent} if it is 
quasi-comapct and quasi-separated. 
We say that a subset $Y$ of a coherent schcme $X$ is a 
{\it Thomason-Ziegler closed subset} if $Y$ is a closed set 
with respect to Zariski topology of $X$ and if $X\ssm Y$ is a 
quasi-compact open subset of $X$. 
The naming of the term `Thomason-Ziegler' comes from the works 
\cite{Zie84} and \cite{Tho97}. 

\begin{para}
\label{df:motives associated with schemes}
{\bf Definition (Motives assocaited with schemes).}\ \ 
Let $X$ be a coherent scheme and $Y$ be a Thomason-Ziegler 
closed subset of $X$. 
We denote the dg-category of perfect complexes 
on $X$ whose cohomological 
supports are in $Y$ by $\Perf_X^Y$. 
In particular we write $\Perf_X$ for $\Perf_X^X$. 
We also write the same letters $(\Perf_X^Y,\qis)$ for 
the Waldhausen category of perfect complexes on $X$ 
whose cohomological supports are in $Y$ 
with the class of quasi-isomorphisms. 
By virtue of \cite[Theorem 0.1.]{Moc13b}, 
the canonical dg-functor 
$\Perf_X^Y\to \calD^b_{\dg}(\Perf_X^Y,\qis)$ is a derived Morita 
equivalence. Thus it induces an isomorphism of motives 
$\cU_{\#}(\Perf_X^Y)\isoto M_{\#}(\Perf_X^Y,\qis)$ for 
$\#\in\{\add,\loc,\nilp \}$. 
We denote it by $M_{\#}^Y(X)$ and call it the {\it additive} 
(resp. {\it localizing, nilpotent invarinat}) {\it motives 
associated with the pair $(X,Y)$}. 
We denote the motives associated with the pair $(X,X)$ by 
$M_{\#}(X)$. 
For a non-negative integer $p$, 
we write $\Perf_X^p$ for the dg-category of perfect complexes on $X$ whose 
cohomological support has codimension$\geq p$ and 
we set $M_{\#,\topo}^p(X):=\cU_{\#}(\Perf^p_X)$ for $\#\in\{\add,\loc,\nilp\}$ 
and call it the {\it topological weight $p$ part of additive} 
(resp. {\it localizing, nilpotent invariant}) {\it motive associated with $X$}.
\end{para}

\begin{para}
\label{rem:fundamental properties of motives}
{\bf Remark (Fundamental properties of motives associated with schemes).}\ \ 
Since the proofs of fundamental properties of algebraic $K$-theory 
in \cite{TT90}, \cite{CHSW08} and \cite{Sch11} are 
based upon derived equivalences of the categories of perfect complexes 
on schemes, the proofs work still fine for motives. 
Namely by replacing $\bbK$ (or $K^B$) with $M_{\#}$ ($\#\in\{\loc,\nilp\}$), 
contiunity (\cite[3.20.2]{TT90}), 
localization (\cite[3.4.9]{Sch11}), excision (\cite[7.1]{TT90}), 
Mayer-Vietoris for Zariski open covers (\cite[3.4.12]{Sch11}) and 
blow up formula for regular center (\cite[3.5.4]{Sch11})
hold for $M_{\#}$.
\end{para}

\begin{para}
\label{df:dash motive}
{\bf Definition (Dash motives associated with noetherian schemes).}\ \ 
Let $X$ be a noetherian scheme and $Y$ be a closed subset of $X$. 
We denote the abelian category of 
coherent sheaves on $X$ whose supports are in $Y$ 
by $\Coh_X^Y$. 
In particular we write $\Coh_X$ for $\Coh_X^X$. 
We set ${M_{\#}'}^Y(X):=\cU_{\#}(\Coh_X^Y)$ for $\#\in\{\add,\loc,\nilp \}$ 
and call it the {\it additive} (resp. {\it localizing, nilpotent invarinat}) 
{\it dash motive associated with the pair $(X,Y)$}. 
We denote the dash motives associated with the pair $(X,X)$ by 
$M_{\#}'(X)$. 
Similarly for an non-negative integer $p$, 
we denote the category of coherent sheaves on $X$ whose support 
has codimension $\geq p$ in $X$ by $\Coh^p_X$. 
We set ${M'}_{\#,\topo}^p(X):=\cU_{\#}(\Coh_X^p)$ for $\#\in\{\add,\loc,\nilp\}$ 
and call it 
{\it topological weight $p$ part of additive} 
(resp. {\it localizing, nilpotent invariant}) 
{\it dash motive associated with $X$}.
\end{para}

\begin{para}
\label{rem:fundamental properties of dash motives}
{\bf Remark (Fundamental properties of dash motives associated with noetherian schemes).}\ \ 
As in Remark~\ref{rem:fundamental properties of motives}, 
since the proofs of fundamental properties of $G$-theory of 
noetherian schemes 
in \cite{TT90} and \cite{Sch11} are 
based upon derived equivalences of 
the categories of pseudo-coherent complexes 
or bounded complexes of coherent modules 
on schemes, the proofs work still fine for dash motives. 
Namely by replacing $G$ (or $K'$) with $M'_{\#}$ ($\#\in\{\loc,\nilp\}$), 
continuity (\cite[3.20.2, 7.2.]{TT90}), 
localization (\cite[3.3.2.]{Sch11}) and excision (\cite[3.19.]{TT90}) 
hold for $M'_{\#}$.
\end{para}

\begin{para}
\label{rem:nilpotent invariance}
{\bf Remark (Nilpotent invariance).}\ \ 
For a coherent scheme $X$, in general $M_{\nilp}(X)$ is not 
isomorphic to $M_{\nilp}'(X_{\red})$. 
But for a noetherian scheme $X$, since the 
inclusion functor $\Coh_{X_{\red}}\rinc\Coh_X$ is 
a nilpotent immersion, it induces an isomorphism 
$M'_{\nilp}(X_{\red})\isoto M'_{\nilp}(X)$. 
The naming of `nilpotent invariant' comes from this fact. 
\end{para}

Specific features of nilpotent invariant motives and 
relationship between motives and dash motives associated with schemes 
summed up with the following proposition. 

\begin{para}
\label{prop:fundamental properties of nilpotent invariant motives}
{\bf Proposition.}\ \ 
{\it Let $X$ be a regular noetherian separated scheme over $\Spec \bbZ$ and $Y$ be a closed subscheme 
of $X$ and $p$ be a non-negative interger. Then
\begin{enumerate}
\enumidef
\item
{\bf (Purity).}\ \ 
The inclusion functor 
$\Coh_Y\rinc\Perf_X^Y$ induces an isomorphism of motives 
$M'_{\nilp}(Y)\isoto M^Y_{\nilp}(X)$. 
In particular if $Y\rinc X$ is an regular closed embedding, then 
we have the caonical isomorphism $M_{\nilp}(Y)\isoto M_{\nilp}^Y(X)$. 
\item
{\bf (Comparison of topological weight).}\ \ 
The inculsion functor 
$\Coh^p_X\rinc \Perf_X^p$ 
induces an isomorphism of 
motives ${M'}_{\nilp,\topo}^p(X)\isoto M_{\nilp,\topo}^p(X)$.
\end{enumerate}
}
\end{para}

\begin{proof}
$\mathrm{(1)}$ 
In $\Mot_{\dg}^{\nilp}$, there exists a commutative diagram of 
localization distinguished triangles 
$$
\xymatrix{
M'_{\nilp}(Y) \ar[r] \ar[d] & M'_{\nilp}(X) \ar[r] \ar[d]_{\textbf{I}} & 
M'_{\nilp}(X\ssm Y) \ar[r] \ar[d]^{\textbf{II}} & \Sigma M'_{\nilp}(Y) \ar[d] \\
M_{\nilp}^Y(X) \ar[r] & M_{\nilp}(X) \ar[r] & M_{\nilp}(X\ssm Y) 
\ar[r] & \Sigma M_{\nilp}^Y(X).
}
$$
Then the morphisms \textbf{I} and \textbf{II} are isomorphism 
by (the proof of) \cite[3.21]{TT90}. 
Thus by five lemma of distinguished triangles, the morphism 
$M'_{\nilp}(Y) \to M_{\nilp}^Y(X)$ is also an isomorphism.  

\sn
$\mathrm{(2)}$ 
Consider the following commutative diagram
$$
\xymatrix{
{M'}^p_{\nilp,\topo}(X) \ar[r] \ar[d] & M^p_{\nilp,\topo}(X) \ar[d]\\
\underset{\codim_X Y\geq p}{\hocolim}M'_{\nilp}(Y)  \ar[r] & 
\underset{\codim_X Y\geq p}{\hocolim} M^Y_{\nilp}(X) 
}
$$
in $\Ho(\Mot_{\dg}^{\nilp})$. 
By contiunity, the vertical morphisms above are isomorphisms 
and by $\mathrm{(1)}$, the bottom morphism is an isomorphism. 
Thus we obtain the result.
\end{proof}

In section~\ref{sec:what makes motives A1homotopy invariance}, 
we will show $\bbA^1$-homotopy invariance and 
projective bundle formula for dash motives of noetherian schemes. 
(See Corollary~\ref{cor:A1-homotopy invariance}.) 

\begin{para}
\label{rem:relative case}
{\bf Remark (Relative version of nilpotent invariant motives).}\ \ 
For a commutative ring $B$ with $1$, 
by replacing $\dgCat$ with $\dgCat_B$ the category 
of small dg-categories over $B$ in the construction 
$\Mot^{\nilp}_{\dg}$, we can obtain 
the stable model category 
$\Mot^{\nilp}_{\dg,B}$ 
of nilpotent invariant motives 
over $B$. 
Similar statements 
as in 
Proposition~\ref{prop:fundamental properties of nilpotent invariant motives} 
are also true for schemes over $\Spec B$. 
\end{para}

\section{Absolute geometric presentation theorem}
\label{sec:absolute geometric presentation theorem}

Gersten's conjecture \cite{Ger73} for Grothendieck groups is equivalent to 
the following generator conjecture. (See \cite{Lev85}, \cite{Dut93}, \cite{Dut95}, \cite{Moc13a} and \cite{Moc16b}.)

\begin{para}
\label{con:gen conj}
{\bf Conjecture (Generator conjecture).}\ \ 
{\it 
For any commutative regular local ring $R$ and 
any natural number $0\leq p \leq \dim R$, 
the Grothendieck group 
$K_0(\Coh^p_{\Spec R})$ is generated by 
cyclic modules $R/(f_1,\cdots,f_p)$ where 
the sequence 
$f_1,\cdots,f_p$ forms an $R$-regular sequence. 
}
\end{para}

For historical background of generator conjecture, 
please see Introduction of \cite{Moc13a}. 
Basically proofs of known cases are based upon structure theorems of $R$ 
(over base). 
For example, if $R$ is smooth over 
a field or a discrete valuation ring, we will utilize a version of 
noether normalization theorem. (See \cite{Qui73} and \cite{GL87}.) 
A version of noether normalization theorem sometimes called 
`geometric presentation theorem'. 
(See for example \cite{C-THK97}.) 
So we would like to formulate an absolute version of such a 
presentation theorem.

On the other hands, by virtue of \cite{GS87} and \cite{TT90}, 
we have the following proposition~\ref{prop:Adams operation}. 
Recall that we say that a scheme $X$ is {\it divisorial} 
if it is quasi-compact and 
if it has an {\it ample family of line bundles}. 
That is, there exists a family of line bundles 
$\{\calL_{\lambda}\}_{\lambda\in\Lambda}$ on $X$ indexed by a non-empty set 
$\Lambda$ which satisfies the following condition. 
For any $f\in\Gamma(X,\calL_{\lambda}^{\otimes m})$, we 
set $X_f:=\{x\in X;f(x)\neq 0 \}$. 
Then the family $\{X_f\}$ is a basis of Zariski topology of $X$ 
where $m$ runs over all positive integer, $\calL_\lambda$ runs over 
the family of line bundles and $f$ runs over all global sections of all of $\calL_{\lambda}^{\otimes m}$. (See \cite[2.2.5.]{Ill71}, \cite[2.1.1.]{TT90} and 
\cite[2.12.]{HM10}.)

\begin{para}
\label{prop:Adams operation} 
{\bf Proposition.}\ \ 
{\it
For a divisorial scheme $X$ and non-negative integer $p\leq \dim X$, 
there exists a family of Adams operations 
$\{\psi_k\}_{k\geq 0}$ indexed by the set of positive integers on 
the Grothendieck group $K_0(\calD(\Perf^p_X))$ 
of the triangulated category of perfect complexes on $X$ 
whose homological support has codimension $\geq p$. 
Moreover if $X$ is the spectrum $\Spec R$ of a noetherian ring $R$ and 
if the sequence $f_1,\cdots,f_p$ is an $R$-regular sequence, then 
we have the equality 
\begin{equation}
\label{eq:Adams operations on Koszul complex}
\psi_k([\kos(f_1,\cdots,f_p)])=k^p[\kos(f_1,\cdots,f_p)]
\end{equation}
for any positive integer $k> 0$ where 
$[\kos(f_1,\cdots,f_p)]$ is the class of $\kos(f_1,\cdots,f_p)$ 
the Koszul complex associated 
with the regular sequencee $f_1,\cdots,f_p$ 
in $K_0(\calD(\Perf_{\Spec R}^p))$. 
}
\end{para}

\begin{proof}
In \cite[4.11]{GS87}, for a positive integer $k$, 
Gille and Soul\'e define an Adams operation of degree $k$
to be any collection of maps 
$\psi_k\colon K_0^Y(X)\to K_0^Y(X)$ 
where $Y$ runs over all closed subsets of $X$ and $K_0^Y(X)$ is 
the Gorthendieck group of relative exact category 
$\Ch_{\geq 0,b}^Y(\cP_X)$ the category of 
bounded complexes $x$ of vector bundles on $X$ 
such that $x_i=0$ for $i<0$ and whose homological support is in $Y$ 
with the class of all quasi-isomorphisms. 
By additivity theorem, 
the degree shift functor 
$$[1]\colon\Ch_{\geq 0,b}^Y(\cP_X)\to 
\Ch_{\geq 0,b}^Y(\cP_X)$$
induces the multiplication by $-1$ on 
the Grothendieck group $K_0^Y(X)$. 
In particular $K_0^Y(X)$ is canonically isomorphic to 
the Grothendieck group of relative exact category of 
the strict perfect complexes(= bounded 
complexes of vector bundles) on $X$ whose homological support 
is in $Y$ with the class of all quasi-isomorphisms. 
Since $X$ is divisorial, this Grothendieck group is isomorphic to 
$K_0(\Perf_X^Y,\qis)$ the Grothendieck group of 
relative exact category of 
perfect complexes on $X$ 
whose homological support is in $Y$ with the class of all quasi-isomorphisms. 
(See \cite[3.8]{TT90}.) 
Since the Adams operations on $K_0^Y(X)$ are compatible with 
the induced map $K_0^Z(X)\to K_0^Y(X)$ 
from the inclusion of closed subsets $Y\rinc Z$ \cite[4.11 A3)]{GS87}, 
it induces the family of Adams operations on 
$$\displaystyle{K_0(\calD (\Perf^p_X))=
\underset{\codim Y\geq p}{\colim} K_0^Y(X)}.$$ 
Finally the equation $\mathrm{(\ref{eq:Adams operations on Koszul complex})}$ 
follows from \cite[4.11 A4')]{GS87}. 
\end{proof}

For a commutative regular ring $R$ and a non-negative integer 
$p\leq \dim R$, 
the inclusion functor $\Coh_{\Spec R}^p\rinc \Perf_{\Spec R}^p$ 
induces an isomorphism 
of the Grothendieck groups 
$K_0(\Coh_{\Spec R}^p)\isoto K_0(\calD(\Perf_{\Spec R}^p))$. 
(Compare Proposition~\ref{prop:fundamental properties of nilpotent invariant motives} $\mathrm{(2)}$ 
and see also 
\cite[Proposition 5.8.]{Moc13a} and \cite[Theorem 3.3.]{HM10}) and via 
this isomorphism, a class $[R/(f_1,\cdots,f_p)]$ of a cyclic module with 
an $R$-sequence $f_1,\cdots,f_p$ in $K_0(\Coh_{\Spec R}^p)$ 
corresponds to the class $[\kos(f_1,\cdots,f_p)]$ of Koszul complex 
$\kos(f_1,\cdots,f_p)$ 
in $K_0(\calD(\Perf^p_{\Spec R}))$. 
Thus roughly saying, 
the generator conjecture says that if $X$ is the spectrum 
$\Spec R$ of a regular local ring $R$, then 
$K_0(\Coh_{\Spec R}^p)$ is generated by objects of Adams weight $p$. 
The purpose of this section is to provide a motivic version of 
the generator conjecture 
and in my viewpoint, it can be regarded as 
an absolute version of a structure theorem of $X$ 
and I would like to call it an absolute geometric presentation 
theorem~\ref{thm:abs geom presentation theorem}. 

We assume that $R$ is regular. Then since 
in particular $R$ is Cohen-Macaulay, 
the ordered set of all ideals of $R$ 
that contains an $R$-regular sequence of length $p$ 
with usual inclusion is directed. 
Thus $\Perf_{\Spec R}^p$ is the filtered limit 
$\underset{I}{\operatorname{colim}}\Perf_{\Spec R}^{V(I)}$ 
where $I$ runs through any 
ideal generated by any $R$-regular sequence of length $p$. 
Thus by continuity of motives, we obtain the following isomorphism
\begin{equation}
\label{eq:continuity}
M_{\#}(\Perf^p_{\Spec R})\isoto \underset{\substack{\codim_{\Spec R}V(I)= p \\
\Spec R/I \underset{\text{regular}}{\hookrightarrow}
\Spec R}}{\hocolim} M_{\#}^{V(I)}(\Spec R)
\end{equation}
for $\#\in\{\loc,\nilp \}$. 
Thus we wish to show that 
for each ideal $I$ generated by $R$-regular sequence 
$f_1,\cdots,f_p$,  
$M_{\nilp}^{V(I)}(\Spec R)$ is spanned by 
objects of Adams weight $p$ in some sense. 

To make the statement more precisely, 
we start by recalling the definitions of Koszul cubes 
from \cite{Moc13a}, \cite{Moc13b} and \cite{Moc16a}. 
Let $S$ be a finite set and $\cC$ be a category. 
An {\it $S$-cube} in $\cC$ is a contravariant functor 
from $\cP(S)$ the power set of $S$ with the usual inclusion order 
to $\cC$. 
A {\it morphism of $S$-cubes} is just a natural transformation. 
We denote the category of $S$-cubes by $\Cub^S(\cC)$. 
For an $S$-cube $x$ and a subset $T$ and an element $t\in T$, 
we write $x_T$ for $x(T)$ and call it a {\it vertex of $x$} ({\it at $T$}) 
and we denote 
$x(T\ssm\{t\}\rinc T)$ by $d^{t,x}_T=d^t_T$ and call it 
a ({\it $k$-direction}) {\it boundary morphism of $x$}. 

Let $S$ be a non-empty finite set such that $\# S=n$ and 
$x$ an $S$-cube in an additive category $\cB$. 
Let us fix a bijection $\alpha$ from $S$ to $(n]$ 
the set of all positive integers $1\leq k\leq n$ 
and we will identify $S$ with the set $(n]$ via $\alpha$. 
We associate an $S$-cube $x$ with 
{\it total complex} $\Tot_{\alpha} x=\Tot x$ 
as follows. 
$\Tot x$ is a chain complex in $\cB$ concentrated in 
degrees $0,\ldots,n$ whose component at degree $k$ is given by 
$\displaystyle{{(\Tot x)}_k:=\underset{\substack{T\in\cP(S) \\ 
\# T=k}}{\bigoplus} x_T}$ 
and whose boundary morphism 
$d_k^{\Tot x}\colon{(\Tot x)}_k \to {(\Tot x)}_{k-1}$ 
are defined by 
$(-1)^{\overset{n}{\underset{t=j+1}{\sum}}\chi_T(t)}d_T^j\colon x_T 
\to x_{T\ssm\{j\}}$ 
on its $x_T$ component to $x_{T\ssm\{j\}}$ component. 
Here $\chi_T$ is the characteristic function $\chi_T\colon S\to \{0,1\}$ of $T$. 
Namely $\chi_T(s)=1$ if $s$ is in $T$ and otherwise $\chi_T(s)=0$. 

Let $A$ be a commutative noetherian ring with $1$. 
By an {\it $A$-sequence} we mean an $A$-regular sequence $f_1,\cdots,f_q$ 
such that any permutation of the $f_j$ is also an $A$-regular sequence. 
A {\it Koszul cube} ({\it associated with an $A$-sequence 
$\ff_S=\{f_s\}_{s\in S}$}) is an $S$-cube in the category of 
finitely generated projective $A$-modules such that for each subset 
$T$ of $S$ and each element $t$ in $T$, $d_T^t$ is an injection and 
$\coker d_T^t$ is annihilated by $f_t^m$ for some positive integer $m$. 
We denote the category of Koszul cubes associated with $\ff_S$ by 
$\Kos_A^{\ff_S}$. 
Notice that if $S$ is the empty set, then $\Kos_A^{\ff_S}$ is just 
the category of finitely generated projective modules $\cP_A$. 
A representative example is the following. 
Let $r$ be a non-negative integer and 
$\fn_S=\{n_s\}_{s\in S}$ a family of non-negative integers 
indexed by $S$ such that $r\geq n_s$ for each $s$ in $S$. 
We define $\Typ_A(\ff_S;r,\fn_S)$ to be an $S$-cube of 
finitely generated free $A$-modules 
by setting for each element $s$ in $S$ and subsets $U \subset S$ and 
$V\subset S\ssm\{s\}$, 
$\Typ_A(\ff_S;r,\fn_S)_U:=A^{\oplus r}$ and 
$\displaystyle{d_{V\sqcup\{s\}}^{\Typ_A(\ff_S;r,\fn_S),s}:=
\begin{pmatrix}
f_s E_{n_s} & 0\\ 
0 & E_{r-n_s} 
\end{pmatrix}}$ 
where $E_m$ is the $m\times m$ unit matrix. 
We call $\Typ_A(\ff_S;r,\fn_S)$ the {\it typical cube of type $(r,\fn_S)$ 
associated with $\ff_S$}. 
We denote the full subcategory of $\Kos_A^{\ff_S}$ consisting of 
those object $\Typ_A(\ff_S;r,\fn_S)$ for some pair $(r,\fn_S)$ 
by $\Kos_{A,\typ}^{\ff_S}$ and call it the {\it category of 
typical Koszul cubes} ({\it associated with $\ff_S$}). 
Roughly saying, $\Kos_{A,\typ}^{\ff_S}$ is just the category of 
objects of Adams weight $\#S$ supported in $V(\ff_S A)$ where $\ff_S A$ is 
an ideal generated by the family $\ff_S$. 

We say that 
a morphism $f\colon x\to y$ in $\Kos_A^{\ff_S}$ is 
a {\it total quasi-isomorphism} if $\Tot f\colon\Tot x\to \Tot y$ is a quasi-isomorphism. 
We denote the class of all total quasi-isomorphisms in 
$\Kos_A^{\ff_S}$ and $\Kos_{A,\typ}^{\ff_S}$ respectively 
by the same letters $\tq$. 
Then the pairs $(\Kos_{A,\typ}^{\ff_S},\tq)$ and $(\Kos_{A,\typ}^{\ff_S},\tq)$ are relative exact categories such that the inclusion functors 
$\Kos_{A,\typ}^{\ff_S}\rinc\Kos_A^{\ff_S}\rinc \Cub^S(\cP_A)$ are 
exact and reflect exactness 
where the last category is the exact category of $S$-cubes in 
the exact category $\cP_A$ of finitely generated projective $A$-modules. 
We also denote the class of all isomorphisms in $\Kos_{A}^{\ff_S}$ 
and $\Kos_{A,\typ}^{\ff_S}$ respectively 
by the same letter $i$. 
The absolute geometric presentation theorem is the following.

\begin{para}
\label{thm:abs geom presentation theorem}
{\bf Theorem (Absolute geometric presentation theorem).}\ \ 
{\it
Let $A$ be a commuative regular ring and $\ff_S=\{f_s\}_{s\in S}$ a 
regular sequence of $A$. Assume that the following two conditions:
\begin{enumerate}
\item
$\ff_S$ is contained in the Jacobson radical of $A$.
\item
For a subset $T$ of $S$, every finitely generated projective 
$A/\ff_TA$-modules are free. 
\end{enumerate}
Then 
\begin{enumerate}
\enumidef
\item
The functor $\Tot\colon (\Kos_{A,\typ}^{\ff_S},\tq)\to 
(\Perf_{\Spec A}^{V(\ff_S)},\qis)$ induces an isomorphism of nilpotent invariant 
motives 
$$M_{\nilp}(\Kos_{A,\typ}^{\ff_S},\tq)\isoto M_{\nilp}^{V(\ff_S)}(\Spec A).$$

\item
The functor $\Tot\colon (\Kos_{A,\typ}^{\ff_S},i)\to 
(\Perf_{\Spec A}^{V(\ff_S)},\qis) $ induces a split epimorphism of 
nilpotent invarinat motives 
$$M_{\nilp}(\Kos_{A,\typ}^{\ff_S})\rdef M_{\nilp}^{V(\ff_S)}(\Spec A).$$ 
\end{enumerate}
}
\end{para}

If $A$ is local, then assumptions in the statement automatically hold. 
If we replace $(\Kos_{A,\typ}^{\ff_S},w)$ with $(\Kos_{A}^{\ff_S},w)$ for $w=\tq$ or $w=i$, then the statement above is essentially proven in \cite[Theorem 0.5, Corollary 9.6.]{Moc13b} and 
\cite[Corollary 5.14.]{Moc13a}. 
(In the references, the statements are 
written in terms of  
$K$-theory. But the same proofs work fine.)
Thus we shall only compare the motives of 
$(\Kos_{A,\typ}^{\ff_S},w )$ and $(\Kos_{A}^{\ff_S},w)$ for $w=\tq$ or $w=i$. 
To give a detailed proof, 
we need intermediate exact categories between $\Kos_{A,\typ}^{\ff_S}$ and $\Kos_{A}^{\ff_S}$. 
We recall some notations from \cite{Moc13a} and \cite{Moc16a}. 
Fix an $S$-cube $x$ 
in an abelian category $\cA$. 
We say that $x$ is {\it monic} if for any pair of subsets 
$U\subset T$ in $S$, $x(U\subset V)$ is a monomorphism. 
For any element $k$ in $S$, 
we define $\Homo_0^k(x)$ to be an $S\ssm\{k\}$-cube in $\cA$ by
setting $\Homo_0^k(x)_T:=\coker d^k_{T\sqcup\{k\}}$ 
for any $T\in\cP(S)$. 
we call $\Homo_0^k(x)$ 
the {\it $k$-direction $0$-th homology} of $x$. 

When $\# S=1$, 
we say that $x$ is {\it admissible} if $x$ is monic, 
namely if its unique boundary morphism is a monomorphism. 
For $\# S>1$, 
we define the notion of an admissible cube inductively 
by saying that 
$x$ is {\it admissible} if $x$ is monic and if 
for every $k$ in $S$, 
$\Homo^k_0(x)$ is admissible. 
For an admissible $S$-cube $x$ and a subset $T=\{i_1,\ldots,i_k\}\subset S$, 
we set 
$\Homo^T_0(x):=\Homo^{i_1}_0(\Homo^{i_2}_0(\cdots(\Homo^{i_k}_0(x))\cdots))$ 
and $\Homo^{\emptyset}_0(x)=x$. 
Notice that $\Homo^T_0(x)$ is an $S\ssm T$-cube for any $T\in\cP(S)$ 
and we can show that 
the definition of $\Homo^T_0(x)$ does not depend upon an order $i_1,\cdots,i_k$ 
up to the canonical isomorphism. 
Let $\fF=\{\cF_T\}_{T\in\cP(S)}$ be 
a family of full subcategories of $\cA$ which are closed under 
isomorphisms. 
Namley for any subset $T$ of $S$ and for any object $x$ in $\cA$ 
such that it is isomorphic to an object in $\cF_T$, $x$ is also an 
object in $\cF_T$. 
We define 
$\ltimes \fF=\underset{T\in\cP(S)}{\ltimes} \cF_T$ 
the {\it multi semi-direct products 
of the family $\fF$} 
as follows. 
$\ltimes \fF$ 
is the full subcategory of $\Cub^S(\cA)$ 
consisting of 
those $S$-cubes $x$ 
such that $x$ is admissible and each vertex of 
$\Homo^T_0(x)$ is in $\cF_T$ for any $T\in\cP(S)$. 

Let $R$ be a commutative noetherian ring with $1$. 
We denote the category of 
finitely generated $R$-modules 
by $\cM_R$.
Let the letter $r$ be a natural number or $\infty$ and 
$I$ be an ideal of $R$. 
Let $\cM_R^I(r)$ be the category of finitely generated $R$-modules 
$M$ such that $\projdim_RM\leq r$ and 
$\supp M\subset V(I)$. 
We write $\cM_R^I$ for $\cM_R^I(\infty)$. 
Since the category is closed under extensions in $\cM_R$, 
it can be considered to be an exact category in the natural way. 
Notice that if $I$ is the zero ideal of $R$, 
then $\cM_R^I(0)$ is 
just the category of finitely generated projective 
$R$-modules $\cP_R$. 
We call an $R$-module $M$ in $\cM_R^I(r)$ is {\it reduced} 
({\it with resect to an ideal $I$}) 
if $IM=0$ and we write $\cM_{R,\red}^I(r)$ for the full subcategory of 
reduced $R$-modules in $\cM_R^I(r)$. 
For $R$-sequence $\ff_S=\{f_s\}_{s\in S}$ and a subset $T$ of $S$, 
we set $\ff_T=\{f_t\}_{t\in T}$ and 
we denote the ideal in $R$ spanned by $\ff_T$ by $\ff_TR$. 
Then we have the following formula
\begin{equation}
\label{eq:Kos ltimes}
\Kos_R^{\ff_S}=\underset{T\in\cP(S)}{\ltimes}\cM_R^{\ff_TR}(\# T).
\end{equation}
(See \cite[4.20]{Moc13a}.) 
Here by convention, 
we set  $\ff_{\emptyset}R=(0)$ the zero ideal of $R$ 
and $\Kos_R^{\ff_{\emptyset}}=\cP_R$ 
the category of finitely generated projective $R$-modules. 
Imitating the equality above, we define the following categories:
$$\Kos_{R,\red}^{\ff_S}:=\underset{T\in\cP(S)}{\ltimes}\cM_{R,\red}^{\ff_TR}(\# T),$$
$$\Kos_{R,\simp}^{\ff_S}:=\underset{T\in\cP(S)}{\ltimes}\cP_{R/\ff_TR}.$$
We call an $S$-cube in $\Kos_{R,\red}^{\ff_S}$ 
(resp. $\Kos_{R,\simp}^{\ff_S}$) {\it reduced} (resp. {\it simple}) 
{\it Koszul cubes associated with $\ff_S$}. 
We denote the class of all total quasi-isomorphisms and 
the class of all isomorphisms in $\Kos_{R,\red}^{\ff_S}$ and $\Kos_{R,\simp}^{\ff_S}$ by 
the same letters $\tq$ and $i$ respectively. 

\begin{proof}[Proof of absolute geometric presntation theorem~\ref{thm:abs geom presentation theorem}]
Under assumpitions 1 and 2 in the statement, 
the inclusion functor 
$\Kos_{A,\typ}^{\ff_S}\rinc \Kos_{A,\simp}^{\ff_S}$ 
is an equivalence of categories. (See \cite[Proposition 1.2.15.]{Moc16a}.) 
Next we will show that the inclusion functor 
$(\Kos_{A,\simp}^{\ff_S} ,w)\to (\Kos_{A,\red}^{\ff_S},w)$ 
induces an isomorphism of 
motives $M_{\#}(\Kos_{A,\simp}^{\ff_S} ,w)\isoto M_{\#}(\Kos_{A,\red}^{\ff_S},w)$ for $w=\tq$ and $w=i$ 
and $\#\in\{\loc,\nilp \}$. 
For the case of $w=i$, it follows from \cite[Corollary 2.1.4.]{Moc16a}. 
(Although the statement in the reference is written in terms of 
non-connective $K$-theory, the same proof works fine.) 
For the case of $w=\tq$, we consider the following commutative square 
of motives in $\Mot_{\dg}^{\#}$ for $\#\in\{\loc,\nilp\}$ 
$$\xymatrix{
M_{\#}(\Kos_{A,\simp}^{\ff_S},\tq) \ar[r] \ar[d]_{M_{\#}(\Homo_0^S)} & 
M_{\#}(\Kos_{A,\red}^{\ff_S},\tq) \ar[d]^{M_{\#}(\Homo_0^S)}\\
M_{\#}(\cP_{R/\ff_S R}) \ar[r] & M_{\#}(\cM_{R,\red}^{\ff_S}(\# S)) 
}$$
where the horizontal morphisms are induced from iclusion functors 
$$\Kos_{A,\simp}^{\ff_S}\rinc \Kos_{A,\red}^{\ff_S}$$
and $\cP_{R/\ff_S R}\rinc \cM_{R,\red}^{\ff_S}(\# S)$. 
The bottom horizontal morphism is isomorphism by \cite[Lemma 1.2.5.]{Moc16b}. 
To show that the vertical morphisms are isomorphisms 
by utilizing \cite[Theorem 8.19]{Moc13b}, we need to check 
assumptions in Theorem 8.19 in \cite{Moc13b}. 
They are follow from \cite[Corollary 5.13.]{Moc13a} for 
the right vertical morphism and \cite[Lemma 2.1.3.]{Moc16a} for the left vertical morphism. 
Finally we will prove that the inclusion functor 
$(\Kos_{A,\red}^{\ff_S},w)\to (\Kos_A^{\ff_S},w)$ 
induces an isomorphism 
$M_{\nilp}(\Kos_{A,\red}^{\ff_S},w)\isoto M_{\nilp}(\Kos_A^{\ff_S},w)$ 
of nilpotent invariant motives for $w=\tq$ or $w=i$.  
The corresponding statement for $K$-theory is Corollary~6.3. in \cite{Moc13a}. 
The proof is based upon derived equivalence and d\'evissage theorem 
for noetherian abelian categories. 
Thus the same proofs work fine for nilpotent invarinat motives. 
We complete the proof. 
\end{proof}

\section{What makes motives $\mathbb{A}^1$-homotopy invariance?}
\label{sec:what makes motives A1homotopy invariance}

In this section, 
we consider what makes motives $\mathbb{A}^1$-homotopy invariance. 
As in Introduction of \cite{MS13}, the main key idea of to prove homotopy 
property of nilpotent invariant dash motives is, roughly speaking, 
that we recognize an affine space as to be a rudimental projective space
\begin{equation}
\label{eq:rudimental projective space}
\bbA^n=\bbP^n\ssm\{[x_0:\cdots:x_n]\in\bbP^n;x_n=0\}(=\bbP^n\ssm\bbP^{n-1}). 
\end{equation}
This means that for a noetherian scheme $X$, 
in the commutative diagram in $\Ho(\Mot^{\nilp}_{\dg})$ below 
\begin{equation}
\label{eq:key commutative diagram}
\xymatrix{
{M'}_{\nilp}(\bbP^{n-1}_X)\ar[r] \ar[d]_{\textbf{I}} & 
{M'}_{\nilp}(\bbP^n_X) \ar[r] \ar[d] &
M'_{\nilp}(X) \ar[r] \ar[d] & 
\Sigma {M'}_{\nilp}(\bbP^{n-1}_X) \ar[d]^{\Sigma\textbf{I}}\\
{M'}_{\nilp}^{\bbP^{n-1}_X}(\bbP^n_X)\ar[r] & 
M'_{\nilp}(\bbP^n_X) \ar[r] &
M'_{\nilp}(\bbA^n_X) \ar[r] & 
\Sigma {M'}_{\nilp}^{\bbP^{n-1}_X}(\bbP^n_X), 
}
\end{equation}
the bottom line is a localization distinguished triangle. 
Moreover by nilpotent invariance of dash motives, 
the vertical morphisms $\textbf{I}$ and $\Sigma\textbf{I}$ 
are isomorphisms. 
We will show that the top line above is also a distinguished triangle 
by a consequence of calculation of nilpotent invariant dash motives of
projective spaces~\ref{cor:nilp inv motives of projective spaces}. 
Thus by five lemma, we obtain homotopy invariance of dash motives. 
Therefore we would like to say that 
$\bbA^1$-homotopy invariance is a degenerated version of 
projective bundle formula and 
nilpotent invariance makes dash motives of noetherian schemes 
$\bbA^1$-homotopy invariance. 
To justify the argument above, 
the main task of this section is to calculate nilpotent invariant 
dash motives of projective spaces. 
In this calucuation, nilpotent inariance of dash motives is again crutial and 
to make the points clarify, we establish an algebraic geometry 
over abelian categories which contains an algebrai geometry over $\bbF_1$ in 
some sense. 
(See Conventions~\ref{conv:F_1-algebra}.) 
We will define for locally noetherian abelian category $\cA$ 
(For definition of locally noetherian abelian category, 
see after the next paragraph.) 
and a positive integer $n$, we define the abelian catgory 
$\bbP^n_{\cA}$ which we call the $n$th projective space over $\cA$. 
The naming is justified by the fact that if $\cA$ is the category 
of quasi-coheret sheaves on a noetherian scheme $X$, 
then $\bbP^n_{\cA}$ is equivalent to the category of 
quasi-coherent sheaves on $\bbP^n_X$ (See Examples~\ref{ex:projective varieties}.) 

To start by defining dash motives associated with 
locally noetherian abelian categories, 
first we recall the conventions and fundamental 
facts of abelian categories from 
\cite{Gab62}, \cite{Pop73}, 
\cite{Ste75}, \cite{Her97} and \cite{TT90}. 
Let $\cA$ be an abelian category and $x$ be an object of $\cA$. 
For a family $\{x_i\rinf x\}_{i\in I}$ 
of subobjects of $x$ indexed by a non-empty set $I$,  
we denonte the minimum object which contains all $x_i$ 
in $\cP(x)$ the partially ordered set of 
isomorphism classes of subobjects of $x$ 
by $\displaystyle{\sum_{i\in I}x_i}$. 
We say that $x$ is {\it finitely generated} if 
for any family of subobjects $\{x_i\}_{i\in I}$ of $x$ such that 
$\displaystyle{x\isoto\underset{i\in I}{\sum}x_i}$, 
there exists a finite subset $J\subset I$ such that 
$\displaystyle{x\isoto\underset{i\in J}{\sum}x_i}$.   
We say that $x$ is {\it finitely presented} if $x$ is finitely generated and 
for any epimorphism $a\colon y\rdef x$ with $y$ finitely generated, 
$\ker a$ is also finitely generated.  
We say that $x$ is {\it coherent} if $x$ is finitely presented and 
if every finitely generated subobject of $x$ is finitely presented. 
We denote the full subcategory of $\cA$ spanned by coherent objects in $\cA$ 
by $\Coh \cA$. 

We say that an abelian category $\cA$ is {\it Grothendieck} if 
$\cA$ has a generator and 
has all small colimits and if 
all small direct limits in $\cA$ is exact. 
The first condition means that 
there exists an object $u$ in a category $\cA$ such that 
the corepresentable functor $\Hom(u,-)$ from $\cA$ to the category of sets 
associated with $u$ is faithful. We call such an object $u$ 
a {\it generator} of $\cA$. 
The last condition means that for any filtered small category $\cI$, 
the colimit functor $\colim_{\cI}\colon\HOM(\cI,\cA) \to \cA$ 
from the category of $\cI$-diagrams in $\cA$ to $\cA$ 
is exact. 
We say that an abelian category $\cA$ is 
{\it locally noetherian} 
if $\cA$ is Grothendieck and if $\cA$ has 
a family of generators consisting of noetherian objects. 
The last condition means that 
there exists a family of objects $\{u_i\}_{i\in I}$ 
in $\cA$ such that 
for each non-zero morphism $a\colon x\to y$ in $\cA$, 
there exists a morphism $b\colon u_i \to x$ for some $i\in I$ 
such that $ab\neq 0$. 
We call such a family a {\it family of generators} of $\cA$. 
If $\cA$ is locally noetherian, 
then every finitely generated object in $\cA$ is noetherian. 
(\cf \cite[Chapter V \S 4]{Ste75}.) 
In particular, every finitely presented (resp. coherent) object in $\cA$ 
is also noetherian. 
Moreover $\Coh\cA$ is a noetherian category. 

For an essentially small abelian category $\cA$, 
we denote the category of left exact functors from 
$\cA^{\op}$ to the category of abelian groups by $\Lex\cA$. 
The category $\Lex\cA$ is a Grothendieck category and there exists 
the Yoneda embedding functor 
$y=y_{\cA}\colon\cA \to \Lex\cA$ which is exact and reflects exactness. 
(See \cite[A.7.14, A.7.16]{TT90}.) 
If $\cA$ is noetherian, then $\Lex\cA$ is locally noetherian and 
the Yoneda functor induces an equivalence of categories 
$\cA\isoto\Coh\Lex\cA$. 
If $\cA$ is locally noetherian, then the inclusion functor 
$\Coh\cA \rinc \cA$, induces an equivalence of categories 
$\Lex\Coh\cA\isoto\cA$. (See \cite[5.8.8, 5.8.9]{Pop73}.)

For example let $X$ be a scheme. 
Then $\Qcoh_X$ the category of quasi-coherent $\cO_X$-modules 
is Grothendieck category. 
(See for \cite[3.14]{Bra14}.) 
Moreover if we assume that $X$ is noetherian, then 
$\Qcoh_X$ is a locally noetherian abelian category and 
the category $\Coh\Qcoh_X$ is just the category of coherent $\cO_X$-modules 
$\Coh_X$ and $\Lex\Coh_X$ is equivalent to $\Qcoh_X$.

\begin{para}
\label{df:dash motives of lna}
{\bf Definition (Dash motives of locally noetherian abelian categories).}\ \ 
Let $\cA$ be a locally noetherian abelian category. 
We define the {\it additive} (resp. {\it localizing, nilpotent invariant}) {\it dash motive of $\cA$} by setting $M'_{\#}(\cA):=M_{\#}(\Coh \cA)$ for $\#\in\{\add,\loc,\nilp \}$. 
\end{para}

\begin{para}
\label{conv:F_1-algebra}
{\bf Conventions ($\bbF_1$-algebra).}\ \ 
In this article, the letters $\bbF_1$ is just a symbol. 
By the term ({\it graded}) {\it $\bbF_1$-algebra}, 
we mean a graded commutative monoid.

For example, $\bbF_1[t_1,\cdots,t_n]$ 
the $n$-variable polynomial ring over $\bbF_1$ is just 
a commutative monoid $\bbN^n$ with usual componentwise addition. 
If we regard it as an $\bbN$-grading algebra, then 
degree $s$ part of $\bbF_1[t_1,\cdots,t_n]$ is given by just the set 
$${\bbF_1[t_1,\cdots,t_n]}_s:
=\left \{(m_1,\cdots,m_n)\in\bbN^n;\sum_{i=1}^nm_i=s \right \}.$$
In this case, we denote an element $(i_1,\cdots,i_n)$ 
in $\bbN^n$ by $t_1^{i_1}\cdots t_n^{i_n}$. 

For a commutative ring $B$ with $1$, 
by a {\it $B$-category} or 
a {\it $B$-functor}, we mean a category and a functor enriched by the category of $B$-modules respectively. 
Similarly by 
an {\it $\bbF_1$-category}, 
an {\it $\bbF_1$-functor} and {\it $\bbF_1$-module}, 
we mean a usual (locally small) 
category, a functor and a set respectively. 

In the rest of this section, let $B$ be a commutative ring with $1$ 
or the letters $\bbF_1$. 

\end{para}

Our approach to define projective spaces over a locally noetherian abelian 
category is based upon the classical result 
of Serre \cite{Ser55} and its globalization by Grothendieck \cite{DG61} and 
its non-commutative version by Artin-Zhang \cite{AZ94} and 
the reconstruction result by Garkusha-Prest \cite{GP08}. 
To recall Serre's theorem, we introduce the notion of stablization of 
categories by an endofunctor. 
For a category $\cC$ 
and a functor $T\colon\cC\to\cC$, 
we define $\stab_T\cC$ to be a category 
by setting $\Ob\stab_T\cC:=\Ob\cC$ and 
$\Hom_{\stab_T\cC}(x,y):=\underset{m\to\infty}{\lim}\Hom_{\cC}(T^mx,T^my)$ 
for any pair of objects $x$ and $y$. We call $\stab_T\cC$ 
the {\it stabilization of $\cC$} ({\it by the functor $T$}). 
If $\cC$ is an abelian category and $T$ is an exact functor, 
then there exists another description of $\stab_T\cC$. Namely 
let $\Nil_T\cC$ be a full subcategory of $\cC$ consisting of those 
objects $x$ such that $T^m(x)=0$ for sufficiently large positive integer $m$. 
Then $\Nil_T\cC$ is a Serre subcategory of $\cC$ and 
for a certain case, there exists a canonical equivalence of categories 
$\stab_T\cC\isoto\cC/\Nil_T\cC$. 
(See Proposition~\ref{prop:fundamental properties of Gr} $\mathrm{(6)}$.) 

For an $\bbN$-grading commutative ring $A$ with $1$, 
we write $\grCoh_A$ for the category of finitely generated graded 
$A$-modules. 
For an integer $k$, there is {\it Serre's twist functor} 
$(k)\colon\grCoh_A \to \grCoh_A$. 
(For definition of Serre's twist functor, see 
Definition~\ref{df:graded cat}. 
Here the definition of abstract Serre's twist functor will be given.)
Here is Serre's theorem. 

\begin{para}
\label{thm:Serre theorem}
{\bf Theorem (Serre).}\ \ (\cf \cite[p.229 Theorem]{AZ94}.)\ \ 
{\it
Let $k$ be a field and let $A$ be an $\bbN$-grading 
commutative ring of finite type over $k$ such that 
it generated by degree one elements. 
Then a functor 
$$\Coh_{\proj A}\to \stab_{(1)}\grCoh_{k}$$
which sends a coherent sheaf $\cF$ on $\proj A$ to 
$\displaystyle{\underset{m\geq 0}{\bigoplus}\Gamma(\proj A,\cF(n))}$ 
gives an equivalence of categories.
}
\qed
\end{para}

To imitate Serre's theorem for more abstract situation, 
first we will define the notion of abstract graded categories 
associated with $\bbN$-grading $B$-algebras and (abelian) $B$-categories.

\begin{para}
\label{nt;graded ring as categeory}
{\bf Notation (Categorified graded ring).} 
Let $A$ be a $\bbN$-grading $B$-algebra. 
We regard $A$ as a category $A_{\graded}$ in the following way. 
The class of objects is just the set of all non-negative integers $\bbN$. 
For any pair of non-negative integers $n$ and $m$, 
we set 
$$\Hom_{A_{\graded}}(n,m):=
\begin{cases}
A_{m-n} & \text{if $m\geq n$}\\
\emptyset & \text{if $m<n$} .
\end{cases}$$
The composition of morphisms is given by the multiplication of $A$.

Similarly, we define ${\bbF_1}_{\graded}$ to be a category by 
setting $\Ob {\bbF_1}_{\graded}:=\bbN$ the set of all non-negative 
integers and for any pair of non-negative integers $n$ and $m$, 
$$
\Hom_{{\bbF_1}_{\graded}}(n,m):=
\begin{cases}
\{\id\} & \text{if $m= n$}\\
\emptyset & \text{if $m\neq n$} .
\end{cases}
$$
\end{para}

\begin{para}
\label{df:graded cat}
{\bf Definition (Graded categories).}\ \ 
Let $A$ be an $\bbN$-grading commutative $B$-algebra 
and $\cC$ be an $B$-category. 
Then we define 
$\gr_A\cC$ to be a category 
by setting 
$$\gr_A\cC:=\HOM_B(A_{\graded},\cC)$$
the category of $B$- functors from 
$A_{\graded}$ to $\cC$ and natural transformations. 
Thus we regard an object $x$ in $\gr_A\cC$ as a system of 
a family of objects 
$\{x_k\}_{k\geq 0}$ in $\cC$ indexed by non-negative integers and 
for each homogeneous element $f\in A_s$, a family of 
morphisms $\{x_k\to x_{k+s} \}_{k\geq 0}$ 
in $\cC$ indexed by non-negative integers which we write the same letter $f$ 
such that this families of morphisms satisfy the equalities arose from 
the relations among elements in $A$. 

We also define $\gr_{\bbF_1}\cC$ 
to be a category 
by setting $\gr_{\bbF_1}\cC:=\HOM_{\bbF_1}({\bbF_1}_{\graded},\cC)$. 
Namely $\gr_{\bbF_1}\cC$ is just a countable products of 
copies of the category $\cC$.

Next we assume that $\cC$ admits a zero object. 
We fix the specific zero object $0$. 
Then for an integer $k$, we define 
{\it Serre's twist functor} 
$
(k)\colon\gr_A\cC\to \gr_A\cC
$ which sends an object $x$ in $\gr_A\cC$  to $x(k)$ in $\gr_A\cC$ 
which defined as follows. 
$$
{x(k)}_n:=
\begin{cases}
x_{m+k} & \text{if $m\geq -k $}\\
0 & \text{if $m<-k$}.
\end{cases}
$$

For a homogeneous element $f\in A_s$ and an object $x$ in $\gr_A\cC$, 
we define $x\to x(s)$ to be a morphism in $\gr_A\cC$ 
which we denote by the same letter $f\colon x\to x(s)$ by 
setting $f\colon x_m\to {x(s)}_m=x_{m+s}$ for each non-negative 
integer $m$. 

For a finite non-empty family $\ff=\{f_k\}_{1\leq k\leq r}$ 
of homogeneous elements in $A$ 
indexed by integers $1\leq k\leq r$, we denote 
the full subcategory of $\gr_A\cC$ consisting of those 
objects $x$ such that there exists a sufficiently large 
integer $N$ such that the morphism $f_k^N\colon x\to x(\deg f_k^N)$ 
is the zero morphism 
for each $1\leq k\leq r$ 
by $\gr_A^{V(\ff)}\cC$. 
Finally let $\gr_{A,\red}^{V(\ff)}\cC$ be a full subcategory of 
$\gr_A^{V(\ff)}\cC$ consisting of those objects $x$ such that 
$f_k\colon x\to x(\deg f_k)$ is the zero morphism for $1\leq k\leq r$. 
\end{para}

\begin{para}
\label{ex:graded cat}
{\bf Examples (Graded category).}\ \ 
\begin{enumerate}
\enumidef
\item
An object 
$\Ob\gr_{\bbF_1[t_1,\cdots,t_n]}\cC\ni x\colon\bbF_1[t_1,\cdots,t_n]\to\cC$ 
can be representated by the diagram
$$
x_0\begin{matrix}\onto{t_1}\\ \vdots\\ \onto{t_n}\end{matrix}x_1
\begin{matrix}\onto{t_1}\\ \vdots\\ \onto{t_n}\end{matrix}x_2
\begin{matrix}\onto{t_1}\\ \vdots\\ \onto{t_n}\end{matrix}\cdots.
$$

\item
Let 
$A$ be $\bbN$-grading commutative $B$-algebra and 
$C$ be a commutative $B$-algebra and let 
$\Mod_C$ be the category of $C$-modules. 
Then the category $\gr_A\Mod_C$ is equivalent to the category 
$\grMod_{A\otimes_B C}$ the category of graded $A\otimes_B C$.

\item
Let $(X,\cO_X)$ be a ringed space and $n$ a positive integer and 
let $\Mod_X$ be the category of $\cO_X$-modules. 
Then the category $\gr_{\bbF_1[t_1,\cdots,t_n]}\Mod_X$ is equivalent to 
the category of graded $\cO_X[t_1,\cdots,t_n]$-modules 
$$\grMod_{\cO_X[t_1,\cdots,t_n]}.$$
\end{enumerate}
\end{para}

Fundamental properties of $\gr_{A}\cA$ is summed up with the following proposition.

\begin{para}
\label{prop:fundamental properties of Gr}
{\bf Proposition (Fundamental properties of $\gr_{A}\cA$).}\ \  
{\it 
Let $\cA$ be a category and $A$ an $\bbN$-grading 
commutative $B$-algebra. Then
\begin{enumerate}
\enumidef
\item
We can calculate a {\rm{(}}co{\rm{)}}limit 
in $\gr_{A}\cA$ by 
term-wise {\rm{(}}co{\rm{)}}limit in $\cA$. 
In particular, if 
$\cA$ is additive {\rm(}resp. abelian, 
closed under small filtered colimits{\rm)} 
then $\gr_A\cA$ is also additive {\rm(}resp. abelian, closed under 
small filtered colimits{\rm)}.

\mn
Moreover assume that $\cA$ is abelian and let $\ff=\{f_k\}_{1\leq k\leq r}$ 
be a family of homogeneous elements of $A$. 
We write $\ff A$ for the ideal of $A$ spanned by $\ff$. Then

\item
$\gr_A^{V(\ff)}\cA$ is a Serre subcategory of $\gr_A\cA$. 

\item
$\gr_{A,\red}^{V(\ff)}\cA$ is a topologizing subcategory 
of $\gr_A^{V(\ff)}\cA$ and the inclusion functor 
$\gr_{A,\red}^{V(\ff)}\cA\rinc \gr_A^{V(\ff)}\cA$ is 
a nilpotent immersion. 

\item
The canonical map $A\rdef A/\ff A$ induces an equivalence of categories 
$\gr_{A/\ff A}\cA\isoto \gr_{A,\red}^{V(\ff)}\cA$. 

\item
The fully faithful embedding $\gr_{A/\ff A}\cA\rinc\gr_A\cA$ 
induced from the canonical map 
$A\rdef A/\ff A$ admits a left adjoint functor
$$-\otimes_A A/\ff A\colon \gr_A\cA \to \gr_{A/\ff A}\cA.$$

\item
The canonical functor 
$\gr_{A}\cA \to
\stab_{(1)}\gr_{A}\cA$ 
induces an isomorphism of categories 
$$\gr_{A}\cA/\Nil_{(1)}\gr_{A}\cA
\isoto \stab_{(1)}\gr_{A}\cA.$$
In particular, $\stab_{(1)}\gr_{A}\cA$ is also an abelian category.

\mn
Moreover assume that $A=B[t_1,\cdots,t_n]/(g_1,\cdots,g_t)$ 
where $B[t_1,\cdots,t_n]$ is the $n$th polynomial ring over $B$ and 
$g_i$ is a homogeneous polynomial in $B[t_1,\cdots,t_n]$ 
for $1\leq i\leq t$. Then

\item
The forgetful functor $U_A\colon \gr_A\cA\to \cA$ which 
sends an object $x$ in $\gr_A\cA$ to $x_0$ admits 
a left adjoint functor $-\otimes_B A\colon \cA \to \gr_A\cA$. 

\item
If $x$ in $\cA$ is a noetherian object, then $x\otimes_B A$ is also 
a noetherian object in $\gr_A\cA$. 

\item
For any object $x$ in $\gr_A\cA$ and for any non-negative integer 
$k\geq 0$, there exists a canonical morphism 
$x_k\otimes_B A(-k)\to x$.

\item
If $\cA$ admits a family of generators 
$\{u_{\lambda}\}_{\lambda\in\Lambda}$ indexed by non-empty set $\Lambda$, 
then $\gr_{A}\cA$ has a family of generators 
$\{u_{\lambda}\otimes_B A(-m)\}_{\lambda\in\Lambda,\ m\geq 0}$. 
In particular if $\cA$ is Grothendieck 
{\rm (}resp. locally noetherian{\rm )}, 
then $\gr_{A}\cA$ is also.

\mn
Moreover assume that $\cA$ is a locally noetherian abelian category.

\item
The inclusion functor $\gr_A\Coh\cA\rinc 
\gr_{A}\cA$ induces an isomorphism of 
categories $\Coh\gr_{A}\Coh\cA\isoto\Coh\gr_{A}\cA$.

\item
Let $x$ be an object in $\gr_A\Coh\cA$. Then 
$x$ is a noetherian object in $\gr_A\Coh\cA$ if and only if 
there exists a positive integer $m>0$ such that 
a canonical morphism 
$\displaystyle{\bigoplus_{k=1}^mx_k\otimes_BA(-k) \rdef x}$ 
induced from the morphisms 
in $\mathrm{(9)}$ is an epimorphism.

\end{enumerate}
}
\end{para}

\begin{proof}
Assertions $\mathrm{(1)}$, $\mathrm{(2)}$, $\mathrm{(3)}$ and 
$\mathrm{(4)}$ are straightfoward. 

\sn
$\mathrm{(5)}$ 
We define $\ff x$ to be a subobject of $x$ by setting 
$$\displaystyle{\ff x:=\sum_{k=1}^r \im (f_k\colon x(-\deg f_k)\to x) }$$
and we set $x\otimes_ A A/\ff A:=x/\ff x$. 
Then we can show that the association which sends an object $x$ in $\gr_A\cA$ 
to an object $x\otimes_A A/\ff A$ gives a left adjoint functor of 
the fully faithful embedding $\gr_{A/\ff A}\cA \rinc \gr_A\cA$. 

\sn
$\mathrm{(6)}$ 
We write $\cB$ and $\cC$ for 
$\gr_{A}\cA/\Nil_{(1)}\gr_{A}\cA$ and 
$\Nil_{(1)}\gr_{A}\cA$ respectively. 
For any pair of objects $x$ and $y$ in $\cB$, we define 
$I_{(x,y)}$ to be a partially ordered direct set by 
$$I_{(x,y)}:=\{(x',y')\in\cP(x)\times\cP(y);x/x',\ y'\in\Ob\cC \}$$
where $\cP(z)$ is the partially ordered set of 
isomorphisms classes of subobjects of $z$ and 
the ordering $(x',y')\leq (x'',y'')$ 
holds if and only if $x''\subset x'$ and $y'\subset y''$. 
Recall that the Hom set $\Hom_{\cB}(x,y)$ is given by the formula
$$\Hom_{\cB}(x,y):=
\underset{(x',y')\in I_{(x,y)}}{\colim}
\Hom_{\gr_{A}\cA}(x',y/y').$$
Let $(x',y')$ be an element in $I_{(x,y)}$. 
Then since $x/x'$ and $y'$ are $(1)$-nilpotent, 
there exists a non-negative integer $m$ such that 
$x/x'(m)\isoto 0$ and $y'(m)\isoto 0$. 
Then we have $(x',y')\leq((x(m))(-m),y')$ and the equalities
\begin{align*}
\Hom_{\gr_{A}\cA}((x(m))(-m),y/y') & =  
\Hom_{\gr_{A}\cA}((x(m))(-m),(y(m))(-m))\\
& =  \Hom_{\gr_{A}\cA}(x(m),y(m)).
\end{align*}
Therefore the canonical map 
$\Hom_{\cB}(x,y)\to\Hom_{\stab_{(1)}\gr_{A}\cA}(x,y)$ 
is an isomorphism. 

\sn
$\mathrm{(7)}$ 
Since $U_A\colon \gr_A\cA\to \cA$ factors through 
$$\gr_A\cA \rinc \gr_{B[t_1,\cdots,t_n]}\cA 
\onto{U_{B[t_1,\cdots,t_n]}} \cA,$$
we shall only prove the case where $A=B[t_1,\cdots,t_n]$ by $\mathrm{(5)}$. 
In this case, for simplicity, for an object $x$ in $\cA$, we denote $x\otimes_B B[t_1,\cdots,t_n]$ 
by $x[t_1,\cdots,t_n]$ which defined as follows. 
For each non-negative integer $s$, degree $s$ part of 
$x[t_1,\cdots,x_n]$ is given by the formula 
$$
{x[t_1,\cdots,t_n]}_s:=\underset{\sum_{i=1}^nm_i=s}{\bigoplus}xt_1^{m_1}\cdots t_n^{m_n}
$$
where $xt_1^{m_1}\cdots t_n^{m_n}$ is just a copy of $x$ and 
for each $1\leq i\leq n$, 
$
x[t_1,\cdots,t_n]_s\to x[t_1,\cdots,t_n]_{s+1}
$ the multiplication by $t_i$ is induced from 
the identity morphisms 
$
\id_x\colon xt_1^{m_1}\cdots t_i^{m_i}\cdots t_n^{m_n}\to xt_1^{m_1}\cdots 
t_i^{m_i+1}\cdots t_n^{m_n}
$ 

Let $k$ be a non-negative integer and $x$ an object 
in $\gr_{B[t_1,\cdots,t_n]}\cA$. 
we define 
$x_k[t_1,\cdots,t_n](-k) \to x$ 
to be a morphism by setting 
for any $m\geq k$ and any $\fraki=(i_1,\cdots,i_n)\in\bbN^n$ 
such that $\displaystyle{\sum_{j=1}^n i_j=m-k}$, 
$$t_1^{i_1}\cdots t_n^{i_n}\colon x_kt_1^{i_1}\cdots t_n^{i_n} \to x_m$$
on the $x_kt^{i_1}\cdots t^{i_n}$ component of 
$x_k[t_1,\cdots,t_n](-k)_m$. 

For any object $x$ in $\cA$ 
and any object $y$ in $\gr_{B[t_1,\cdots,t_n]}\cA$, 
we have a functorial isomorphism 
$\Hom_{\cA}(x,y_0)\isoto
\Hom_{\gr_{B[t_1,\cdots,t_n]}\cA}(x[t_1,\cdots,t_n],y)$ 
which sends $f$ to 
$(x[t_1,\cdots,t_n]\onto{f\otimes_B B[t_1,\cdots,t_n]}y_0[t_1,\cdots,t_n] \to y)$.

\sn
$\mathrm{(8)}$ 
If $A=B[t_1,\cdots,t_n]$, assertion is proven in \cite[4.14 (1)]{MS13}. 
(In the reference, the proof is written for $B=\bbF_1$ with 
slightly different conventions. But a similar proof works fine for a general 
$B$.) 
In a general case, $x\otimes_B A$ is a quotinet of 
a noetherian object $x[t_1,\cdots,x_n]$ by 
$\fg x[t_1,\cdots,t_n]$ where $\fg=\{g_i\}_{1\leq i\leq t}$. 
Thus $x\otimes_B A$ is also a noetherian object. 

\sn
$\mathrm{(9)}$ 
We regard $x$ as an object of $\gr_{B[t_1,\cdots,t_n]}\cA$. 
Then as in the proof of $\mathrm{(7)}$, 
there exists a canonical morphism $x_k[t_1,\cdots,t_n](-k)\rdef x$.
Since $x$ is in $\gr_A\cA$, it factors through 
$x_k[t_1,\cdots,t_n](-k)\rdef x_k\otimes_B A(-k)\to x$ and 
the last morphism is the desired morphism.  

\sn
$\mathrm{(10)}$ 
For any non-zero morphism $a\colon x\to y$ in $\gr_A\cA$, there exists a non-negative integer $k$ such that $a_k\colon x_k\to y_k$ is a non-zero morphism. 
Then there exists $\lambda\in\Lambda$ and 
$b\colon u_{\lambda}\to x_k$ such that 
$a_kb\neq 0$. 
Then the compositions 
$u_\lambda\otimes_B A(-k) \onto{b\otimes_B A(-k)} 
x_k\otimes_B A(-k) \to x 
\onto{a} y$ is also a non-zero morphism. 
Therefore $\{u_{\lambda}\otimes_B A(-m)\}_{\lambda\in\Lambda,\ m\geq 0}$ 
is a family of generators of $\gr_{A}\cA$. 
The last assertion follows from $\mathrm{(1)}$ and $\mathrm{(8)}$. 

\sn
$\mathrm{(11)}$ 
We need to prove that for any object 
$x$ in $\Coh\gr_{A}\cA$, 
$x_k$ is in $\Coh\cA$ for any $k\geq 0$. 
Obviously if $x$ is noetherian, then $x(k)$ is also noetherian. 
Therefore we shall assume that $k=0$. 
Let $y_0\rinc y_1\rinc \cdots \rinc y_n\rinc \cdots$ 
be a sequence of subobjects in $x_0$, 
then we define subobject $\bar{y}_k$ of $x$ as follows:
$${(\bar{y}_k)}_l:=
\begin{cases}
y_k & \text{for $l=0$}\\
x_l & \text{for $l\geq 1$}.
\end{cases}$$
Then since $x$ is noetherian, 
there exists a positive integer $m$ such that 
$\bar{y}_m=\bar{y}_{m+1}=\cdots $ and this means that 
$y_m=y_{m+1}=\cdots$. 
Hence $x_0$ is a noetherian object in $\cA$.

\sn
$\mathrm{(12)}$ 
Assume that $x$ is a noetherian object. 
For each integer $m\geq 0$, 
we set $\displaystyle{y_i:=\im\left (\bigoplus_{k=0}^i x_k\otimes_BA(-k) \to x \right ) }$. 
Then the asscending filtration $\{y_i \}_{i\geq 0}$ is stational. 
Say $y_m=y_{m+1}=\cdots$. 
Since $y_{\infty}=x$, we obtain the equality $y_m=x$ and this means 
that $m$ is the desired integer.  

\sn
Conversely if there exists an epimorphism 
$\displaystyle{\bigoplus_{k=1}^mx_k\otimes_BA(-k) \rdef x}$, then 
since $x$ is a quotinet of finite direct sum of noetherian objects, 
$x$ is also noetherian.

\end{proof}

\begin{para}
\label{thm:motive of grF1}
{\bf Theorem (Dash motive of $\gr_{\bbF_1[t_1,\cdots,t_n]}\cA$).}\ \ 
(\cf \cite[4.24]{MS13}.)\ \ 
Let $\cA$ be a locally noetherian abelian category. 
Then for $\#\in\{\loc,\nilp\}$, we have 
the canonical isomorphism
\begin{equation}
\label{eq:Cohgr isom}
M'_{\#}\left(\gr_{\bbF_1[t_1,\cdots,t_n]}\cA\right)\isoto
\underset{i\geq 0}{\bigoplus}M'_{\#}(\cA)s^i
\end{equation}
where $M'_{\#}(\cA)s^i$ is a copy of $M'_{\#}(\cA)$ and 
we denote the identity morphism 
$\id_{M'_{\#}(\cA)}\colon M'_{\#}(\cA)s^j\to M'_{\#}(\cA)s^{j+i}$
by $\times s^i$ 
(we write $\times s$ for $\times s^1$.) 
Then the isomorphism $\mathrm{(\ref{eq:Cohgr isom})}$
makes the diagram below commutative. 
\begin{equation}
\label{eq:comm square of Cohgr isom}
\xymatrix{
M'_{\#}\left(\gr_{\bbF_1[t_1,\cdots,t_n]}\cA\right) \ar[r]^{\ \ \ \ \ \ \sim} \ar[d]_{(-1)} & 
\underset{i\geq 0}{\bigoplus}M'_{\#}(\cA)s^i \ar[d]^{\times s}\\
M'_{\#}\left(\gr_{\bbF_1[t_1,\cdots,t_n]}\cA\right) \ar[r]^{\ \ \ \ \ \ \sim} & 
\underset{i\geq 0}{\bigoplus}M'_{\#}(\cA)s^i.
}
\end{equation}
\end{para}

To give a proof of Theorem~\ref{thm:motive of grF1}, 
we need to recall several concepts from \cite{MS13} with slightly changing 
notations. 
Let $\cA$ be an abelian category and let 
$x$ be an object in $\gr_{\bbF_1[t_1,\cdots,t_n]}\cA$. 
For each non-negative integer $m\geq 0$, 
we define $F_mx$ to be a subobject of $x$ by setting 
$$
{(F_mx)}_k:=
\begin{cases}
x_k & \text{if $k\leq m$}\\
\displaystyle{\sum_{\sum_{j=1}^ni_j=k-m}\im t_1^{i_1}\cdots t_n^{i_n}}& \text{if $k>m$}.
\end{cases}
$$
By convention, we set $F_{-1}x=0$. 
We call a family $\{F_m x\}_{m\geq 0}$
a {\it canonical filtration of $x$}. 
Moreover if we assume that $x$ is noetherian, 
then there exists a sufficiently large positive integer $m$ 
such that $x=F_mx$ and we set 
$\deg x:=\min\{m\in\bbN;x=F_mx \}$ 
and call it {\it degree of $x$}. 

We set $(n]:=\{k\in\bbN;1\leq k\leq n \}$. 
We define 
$
\kos(x)\colon {\cP((n])}^{\op}\to \gr_{\bbF_1[t_1,\cdots,t_n]}\cA
$ 
to be an $(n]$-cube by sending a subset $T$ in $\cP((n])$ to $x(-\# T)$, 
and an inclusion 
$T\ssm\{s\} \rinc T$ to $x(-\# T) \onto{t_s} x(-\# T+1)$. 
We call $\kos(x)$ the {\it Koszul cube associated with $x$}. 
Moreover we define $T_i(x)$ to be the {\it $i$th Koszul homology of $x$} 
by setting 
$T_i(x):=\Homo_i(\Tot\kos(x))$. 
We say that $x$ is {\it $t$-regular} if $T_i(x)=0$ for any $i>0$. 

\begin{proof}[Proof of Theorem~\ref{thm:motive of grF1}]
By replacing $\cA$ with $\Coh\cA$ and by 
Proposition~\ref{prop:fundamental properties of Gr} $\mathrm{(3)}$, 
we shall assume that $\cA$ is noetherian. 
Let 
${\Coh\gr_{\bbF_1[t_1,\cdots,t_n]}\cA}_{t-\reg}$ 
be a full subcategory of 
$\Coh\gr_{\bbF_1[t_1,\cdots,t_n]}\cA$ 
consisting of $t$-regular objects. 
Then the inclusion functor 
${\Coh\gr_{\bbF_1[t_1,\cdots,t_n]}\cA}_{t-\reg}\rinc
\Coh\gr_{\bbF_1[t_1,\cdots,t_n]}\cA$ 
induces an equivalence of triangulated categories 
$$
\calD^b\left({\Coh\gr_{\bbF_1[t_1,\cdots,t_n]}\cA}_{t-\reg}\right)\isoto
\calD^b\left(\Coh\gr_{\bbF_1[t_1,\cdots,t_n]}\cA\right)
$$
on bounded derived categories. 
(See \cite[4.24 (1)]{MS13}.) 
Next for each non-negative integer $m$, 
let 
${\Coh\gr_{\bbF_1[t_1,\cdots,t_n]}\cA}_{t-\reg,\deg\leq m}$ 
be the full subcategory of 
${\Coh\gr_{\bbF_1[t_1,\cdots,t_n]}\cA}_{t-\reg}$ 
consisting of those objects of degree$\leq m$. 
We define 
$$a\colon \cA^{\times m+1}\to 
{\Coh\gr_{\bbF_1[t_1,\cdots,t_n]}\cA}_{t-\reg,\deg\leq m}\ \ \text{and}$$ 
$$b\colon {\Coh\gr_{\bbF_1[t_1,\cdots,t_n]}\cA}_{t-\reg,\deg\leq m}
\to \cA^{\times m+1}$$
to be the exact functors 
by sending an object ${(x_k)}_{0\leq k\leq m}$ to an object 
$\displaystyle{\bigoplus_{k=0}^{m}x_k[t_1,\cdots,t_n](-k)}$ and an object 
$x$ to an object 
${({T_0(x)}_k)}_{0\leq k\leq m}$ 
respectively. Well-definedness of this functors 
follows from \cite[4.14, 4.23 (1)]{MS13}. 
Then we have equalities
$$
ab(x)=\bigoplus_{k=0}^m{T_0(x)}_k[t_1,\cdots,t_n](-k)\underset{\text{\cite[4.23]{MS13}}}{\isoto}
\bigoplus_{k=0}^mF_px/F_{p-1}x,
$$
$$
ba({(x_k)}_{0\leq k\leq m})={\left ({T_0\left (\bigoplus_{k=0}^m 
x_k[t_1,\cdots,x_n](-k)
\right)}_k \right )}_{0\leq k\leq m}\isoto 
{(x_k)}_{0\leq k\leq m}.
$$
Thus by additivity, in $\Ho(\Mot^{\#}_{\dg})$, 
$ab$ and $ba$ induce the identity morphisms on dash motives. 
On the other hand, since we have the equality 
$\displaystyle{
{\Coh\gr_{\bbF_1[t_1,\cdots,t_n]}\cA}_{t-\reg}
\isoto\underset{m\to\infty}{\colim}\ 
{\Coh\gr_{\bbF_1[t_1,\cdots,t_n]}\cA}_{t-\reg,\deg\leq m}
}$, 
we have isomorphisms in $\Ho(\Mot^{\#}_{\dg})$:
\begin{eqnarray*}
\label{eq:Cohgr isomorphisms}
M'_{\#}\left(\gr_{\bbF_1[t_1,\cdots,t_n]}\cA\right)
& \underset{\text{Morita equivalence}}{\isoto} &
M'_{\#}\left({\gr_{\bbF_1[t_1,\cdots,t_n]}\cA}_{t-\reg}
\right) \nonumber \\
& \underset{\text{continuity}}{\isoto} &
\underset{m\to\infty}{\colim}\ 
M'_{\#}\left({\gr_{\bbF_1[t_1,\cdots,t_n]}\cA}_{t-\reg,\deg\leq m}\right) \nonumber \\
& \underset{\text{additivity}}{\isoto} &
\underset{m\to\infty}{\colim}\ 
\bigoplus_{i=0}^mM'_{\#}(\cA)s^i \nonumber\\
& \isoto & \bigoplus_{i=0}^{\infty}M'_{\#}(\cA)s^i.
\end{eqnarray*}
Commutativity of the diagram $\mathrm{(\ref{eq:comm square of Cohgr isom})}$ 
follows from the definition of the functor $a$ and the following 
commutative diagram.
$$
\xymatrix{
& M'_{\#}(\cA) \ar[ld]_{\times s^i} \ar[rd]^{\ \ \ \ \ -\otimes\bbF_1[t_1,\cdots,t_n](-i)} \\
M'_{\#}(\cA)s^i \ar[rr] & & 
M'_{\#}\left (\gr_{\bbF_1[t_1,\cdots,t_n]}\cA \right ).
}
$$
\end{proof}

In the light of Theorem~\ref{thm:Serre theorem} and 
Proposition~\ref{prop:fundamental properties of Gr} $\mathrm{(4)}$, 
we would like to consider a quotient category of $\gr_A\cA$ by 
$\Nil_{(1)}\gr_A\cA$. 
To examine quotient categories of abelian category, 
we recall the quotient theory of locally coherent abelian categories 
from \cite{Her97}. 
We say that a Grothendieck category $\cA$ is {\it locally coherent} 
if every object of $\cA$ is a direct limit of coherent objects of $\cA$. 
The following lemma is fundamental. 

\begin{para}
\label{lem:loc noetherian}
{\bf Lemma.} 
{\it
A locally noetherian category $\cA$ is locally coherent.
}
\end{para}

\begin{proof}
Let $\{u_i\}_{i\in I}$ be 
a family of noetherian generators of $\cA$ 
and let $x$ be an object of $\cA$. 
For a non-zero morphism $\alpha\colon x\to y_{\alpha}$, 
there exists $i_{\alpha}$ in $I$ and a morphism 
$\beta_{\alpha}\colon u_{i_{\alpha}}\to x$ such that 
$\alpha\beta_{\alpha}\neq 0$. 
Then the canonical morphism 
$\displaystyle{\bigoplus_{\alpha}u_{i_{\alpha}}\to x}$ 
induced from the family $\{\beta_{\alpha}\colon u_{i_{\alpha}}\to x\}_{\alpha}$ 
is epimorphism and we have the canonical isomorphism 
$\displaystyle{x\isoto \sum_{\alpha}\im \beta_{\alpha}}$. 
Here $\im\beta_{\alpha}$ is noetherian, 
a fortiori, coherent for any $\alpha\colon x\to y_{\alpha}$. 
Therefore $\cA$ is locally coherent.
\end{proof}

Let $\cA$ be a Grothendieck category 
and $\cB$ a full subcategory of $\cA$. 
We say that 
$\cB$ is a  {\it localizing subcategory} 
({\it of $\cA$}) if 
$\cA$ is closed under sub- and quotient objects and extensions and coproducts. 
(In \cite[2.1]{Her97}, we call a localizing subcategory a {\it hereditary torsion subcategory}.) 
We write ${}^L\!\!\!\sqrt{\cB}$ 
for intersection of all localizing subcategories 
which contain $\cB$ and call it the 
{\it localizing radical of $\cB$} ({\it in $\cA$}). 
Assume that $\cA$ is locally coherent and $\cB$ is localizing. 
We say that $\cB$ is of finite type if 
there exists a Serre subcategory $\cT$ of $\Coh\cA$ 
such that $\cB={}^L\!\!\!\sqrt{\cT}$.

Let $\cA$ be a locally coherent abelian category. 
There exists an inclusion preserving bijiective correspondence 
between the class of Serre subcategories of $\Coh\cA$ and 
the class of localizing subcategories of $\cA$ of finite type. 
The correspondence given by sending 
a Serre subcategory
$\cB$ of $\Coh\cA$ to 
${}^L\!\!\!\sqrt{\cB}$ 
and 
a localizing subcategory 
$\cT$ of $\cA$ of finite type to 
$\cT\cap\Coh\cA$ 
which are mutual inverse. 
(See \cite[Theorem 2.8.]{Her97}.) 

Let $\cA$ be a locally coherent category and 
$\cB$ a Serre subcategory of $\Coh\cA$. 
The inclusion functor 
${}^L\!\!\!\sqrt{\cB}\rinc\cA$ admits a right adjoint functor 
$t_{{}^L\!\!\!\sqrt{\cB}}\colon\cA \to {}^L\!\!\!\sqrt{\cB}$ and we call it 
the {\it torsion functor} ({\it associated to ${}^L\!\!\!\sqrt{\cB}$}). 
$t_{{}^L\!\!\!\sqrt{\cB}}$ is a left exact functor and 
we write $t_{{}^L\!\!\!\sqrt{\cB}}^1$ for 
the first derived functor of $t_{{}^L\!\!\!\sqrt{\cB}}$. 
We say that an object $x$ of $\cA$ is 
{\it ${}^L\!\!\!\sqrt{\cB}$-torsion free} 
(resp. {\it ${}^L\!\!\!\sqrt{\cB}$-closed}) 
if $t_{{}^L\!\!\!\sqrt{\cB}}(x)=0$. 
(resp. $t_{{}^L\!\!\!\sqrt{\cB}}(x)=0$ and $t_{{}^L\!\!\!\sqrt{\cB}}^1(x)=0$.) 
We denote the full subcategory of $\cA$ 
spanned by ${}^L\!\!\!\sqrt{\cB}$-closed objects 
by $\cA/{}^L\!\!\!\sqrt{\cB}$ and call it the 
{\it quotient category of $\cA$ by ${}^L\!\!\!\sqrt{\cB}$}. 
$\cA/{}^L\!\!\!\sqrt{\cB}$ is a locally coherent category. 
The inclusion functor 
$j_{{}^L\!\!\!\sqrt{B}}\colon\cA/{}^L\!\!\!\sqrt{\cB}\rinc \cA$ 
admits an exact left adjoint functor 
$q_{{}^L\!\!\!\sqrt{\cB}}\colon\cA \to \cA/{}^L\!\!\!\sqrt{\cB}$ 
such that 
$q_{{}^L\!\!\!\sqrt{\cB}}j_{{}^L\!\!\!\sqrt{B}}=\id_{\cA/{}^L\!\!\!\sqrt{B}}$. 
The functor $q_{{}^L\!\!\!\sqrt{\cB}}$ induces an exact functor 
$\Coh q_{{}^L\!\!\!\sqrt{\cB}}\colon\Coh\cA\to
\Coh\left (\cA/{}^L\!\!\!\sqrt{\cB} \right )$ 
and an equivalence of categories 
$\left (\Coh\cA \right )/\cB\isoto\Coh\left (\cA/{}^L\!\!\!\sqrt{\cB} \right )$. 
Here the category $\left (\Coh\cA\right )/\cB$ is 
the usual quotient abelian category of $\Coh\cA$ 
by the Serre subcategory $\cB$. 
(\cf \cite[\S 2.2]{Her97}.) 

\begin{para}
\label{prop:quot of loc noe}
{\bf Proposition.} 
{\it
Let $\cA$ be a locally noetherian category and 
$\cB$ a Serre subcategory of 
$\Coh\cA$. 
Then
\begin{enumerate}
\enumidef
\item
The localization functor 
$q_{{}^L\!\!\!\sqrt{\cB}}\colon \cA\to\cA/{}^L\!\!\!\sqrt{\cB}$ sends 
noetherian objects in $\cA$ to noetherian objects in $\cA/{}^L\!\!\!\sqrt{\cB}$.

\item
The quotient category 
$\cA/{}^L\!\!\!\sqrt{\cB}$ is also locally noetherian.
\end{enumerate}
}
\end{para}

\begin{proof}
For simplicity we write $q$ and $j$ for $q_{{}^L\!\!\!\sqrt{\cB}}$ 
and $j_{{}^L\!\!\!\sqrt{\cB}}$ respectively.

\sn
$\mathrm{(1)}$ 
Let $x$ be a noetherian object in $\cA$. 
Then $q(x)$ is coherent in $\cA/{}^L\!\!\!\sqrt{\cB}$. 
In particular $jq(x)$ 
is finitely generated in $\cA$ 
by \cite[2.2]{Her97}. 
Since $\cA$ is locally noetherian, 
$jq(x)$ 
is $\cA$-noetherian. 
In particular $q(x)$ 
is noetherian in $\cA/{}^L\!\!\!\sqrt{\cB}$. 

\sn
$\mathrm{(2)}$ 
First notice that 
$\cA/{}^L\!\!\!\sqrt{\cB}$ is Grothendieck. 
Let $\{u_{\lambda}\}_{\lambda\in\Lambda}$ be 
a fimily of noetherian generators of $\cA$. 
We will prove that $\{q(u_{\lambda})\}_{\lambda\in\Lambda}$ 
is a family of genrators of $\cA/{}^L\!\!\!\sqrt{\cB}$. 
Then since $q(u_{\lambda})$ is noetherian 
in $\cA/{}^L\!\!\!\sqrt{\cB}$ 
for any $\lambda$ in $\Lambda$ 
by $\mathrm{(1)}$, 
we complete the proof. 
Let $\theta\colon\id_{\cA}\to jq$ be an adjunction morphism and 
$\beta\colon x\to y$ be a non-zero morphism in $\cA/{}^L\!\!\!\sqrt{\cB}$. 
Then there exist an element $\lambda$ in $\Lambda$ and 
a morphism $\alpha\colon u_{\lambda} \to j(x)$ such that 
$j(\beta)\alpha\neq 0$. 
Then by ajointness of $q$, there exists a morphism 
$\alpha'\colon q(u_{\lambda}) \to x$ such that 
$\alpha=\theta(u_{\lambda})j(\alpha')$. 
Obviously we have $\beta\alpha'\neq 0$. 
Hence we obtain the result.
\end{proof}

Now we give a definition of projective varieties over a
locally noetherian abelian categories. 

\begin{para}
\label{df:projective varieties}
{\bf Definition (Projective varieties).} 
Let $A$ be a $\bbN$-grading $B$-algebra and 
let $\cA$ be a locally noetherian abelian $B$-category. 
Then we set 
$\proj_{\cA} A:=\gr_A\cA/{}^L\!\!\!\sqrt{\Nil_{(1)}\gr_A\cA}$ 
and call it {\it proj $A$ over $\cA$}.

Moreover let $\ff=\{f_k\}_{1\leq k\leq r}$ is a family of 
homogeneous elements in $A$. 
Then we set $\proj_{\cA}^{V(\ff)} A:=\gr_A^{V(\ff)}\cA/{}^L\!\!\!\sqrt{\Nil_{(1)}\gr_A\cA}$. 

By virtue of Proposition~\ref{prop:quot of loc noe} $\mathrm{(2)}$ and 
\cite[\S 2]{Her97}, 
$\proj_{\cA} A$ and $\proj_{\cA}^{V(\ff)} A $ 
are locally noetherian and we 
have equalities 
$$\Coh\proj_{\cA} A= 
\Coh \gr_A\cA/\Nil_{(1)}\Coh\gr_A\cA \ \ \text{and}$$ 
$$\Coh\proj_{\cA}^{V(\ff)} A= 
\Coh \gr_A^{V(\ff)}\cA/\Nil_{(1)}\Coh\gr_A^{V(\ff)}\cA.$$ 

Thus $\proj A$ is a functor from the (large) category of 
locally noetherian abelian categories and cocontinuous functors to 
itself.

In particular, for a non-negative integer $n$, 
we set 
$\bbP^n_{\cA}:=\proj_{\cA} \bbF_1[t_0,\cdots,t_n]$ 
and call it the {\it $n$th projective space over $\cA$}. 

For $\#\in\{\add,\loc,\nilp\}$, 
we write ${M'}_{\#}^{V(\ff)}(\proj_A\cA)$ 
for $M(\Coh\proj_{\cA}^{V(\ff)}\cA)$ a dash motive of 
$\proj_{\cA}^{V(\ff)} A$. 
\end{para}

\begin{para}
\label{ex:projective varieties}
{\bf Examples (Projective varieties).}\ \ 
\begin{enumerate}
\enumidef
\item
(\cf \cite[p.49]{Mum66}.) 
Let $X$ be a noetherian scheme and $n$ a positive integer. 
Then we have the canonical equivalence of categories 
$\bbP^n_{\Qcoh_X} \isoto \Qcoh_{\mathbb{P}^n_X}$.

\item
Let $k$ be a field, $A$ an $\mathbb{N}$-graded commutative ring 
of finite type over $k$ such that $A$ is generated by elements of degree one. 
Then by Serre's theorem~\ref{thm:Serre theorem}, 
we have the canonical equivalence of categories 
$\proj_{\Mod_k}A\isoto\Qcoh_{\proj A}$.
\end{enumerate}
\end{para}

To study the localization sequence of 
dash motives associated with 
projective varieties 
over a locally noetherian abelian categroy and calculate 
the nilpotent invariant dash motives associated with projective spaces, 
we recall the notion of special filtering Serre subcategories. 
Let $\cS$ be a Serre subcategory of an abelian category $\cA$. 
We say that $\cS$ is ({\it right}) {\it special filtering} if 
for any monomorphism $x\rinf y$ in $\cA$ with $x$ in $\cS$, 
there exists a morphism $y \to z$ with $z$ in $\cS$ such that 
the composition $x \to y \to z$ is a monomorphism. 
(\cf \cite[Definition 1.5.]{Sch04})
If $\cS$ is special filtering, 
then the canonical sequence 
$$\calD^b_{\dg}(\cS) \to \calD^b_{\dg}(\cA) \to \calD^b_{\dg}(\cA/\cS)$$
is an exact sequence of dg-categories. 
(See \cite[4.1]{Gro77} and \cite[1.15]{Kel99}.) 
Thus it induces a distinguished triangle 
$$M'_{\#}(\cS)\to M'_{\#}(\cA) \to M'_{\#}(\cA/\cS)\to \Sigma M'_{\#}(\cS)$$
in $\Ho(\Mot^{\#}_{\dg})$ for $\#\in\{\loc,\nilp \}$. 

The following lemma is fundamental. 

\begin{para}
\label{lem:Nil-> Gr-> proj}
{\bf Lemma.}\ \ 
{\it
Let $A$ be an $\bbN$-grading $B$-algebra and $\ff=\{f_k\}_{1\leq k\leq r}$ 
be a family of homogeneous elements in $A$ 
and $\cA$ be a locally 
noetherian abelian category. 
Then $\Nil_{(1)}\Coh\gr_A\cA$ and $\Nil_{(1)}\Coh\gr_A^{V(\ff)}\cA$ are 
right special filtering in 
$\Coh\gr_A\cA$ and 
$\gr_A^{V(\ff)}\cA$ respectively. 
In particular we have the distinguished triangules 
$$M'_{\#}(\Nil_{(1)}\gr_A\cA)\to M'_{\#}(\gr_A\cA)\to 
M'_{\#}(\proj_{\cA}A) \to 
\Sigma M'_{\#}(\Nil_{(1)}\gr_A\cA)\ \ \text{and}$$ 
$$M'_{\#}(\Nil_{(1)}\gr_A^{V(\ff)}\cA)\to M'_{\#}(\gr_A^{V(\ff)}\cA)\to 
M'_{\#}(\proj_{\cA}^{V(\ff)}A) \to 
\Sigma M'_{\#}(\Nil_{(1)}\gr_A^{V(\ff)}\cA) $$ 
in $\Mot_{\dg}^{\#}$ for $\#\in\{\loc,\nilp\}$. 
}
\end{para}

\begin{proof}
Let $x\rinf y$ be a monomorphism in 
$\Coh\gr_A\cA$ with $x$ in $\Nil_{(1)}\Coh\gr_A\cA$. 
Then there exists a non-negative integer 
$k\geq 0$ such that $x_m=0$ for $m\geq k$. 
Then we define $z$ to be an object in $\Nil_{(1)}\Coh\gr_A\cA$ 
by setting $z_l=y_l$ if $l< k$ and $z_l=0$ if $l\geq k$. 
There is a canonical epimorphism $y\rdef z$ such that the composition 
$x \rinf y \rdef z$ is a monomorphism. 
Thus $\Nil_{(1)}\Coh\gr_A\cA$ is right special filtring in $\Coh\gr_A\cA$. 
A proof for $\Nil_{(1)}\Coh\gr_A^{V(\ff)}\cA$ is similar. 
\end{proof}

\begin{para}
\label{cor:nilp inv motives of projective spaces}
{\bf Corollary (Nilpotent invariant dash motive of projective spaces).}\ \ 
{\it
Let $\cA$ be 
a locally noetherian abelian category and let $n$ be a non-negative integer. 
Then the exact functor 
$\displaystyle{\bigoplus_{i=0}^n\Coh\cA\to\Coh\bbP^n_{\cA}}$ 
which sends an object ${(x_i)}_{0\leq i\leq n}$ to 
$\displaystyle{\bigoplus_{i=0}^nx_i[t_1,\cdots,t_n](-i)}$ 
induces an isomorphism of 
nilpotent invariant dash motives 
$\displaystyle{M'_{\nilp}\left ( \bbP^n_{\cA}\right)\isoto \bigoplus_{i=0}^n 
M'_{\nilp}(\cA)s^i}$ 
where $M'_{\nilp}(\cA)s^i$ is just a copy of $M'_{\nilp}(\cA)$.
}
\end{para}

The proof is carried out in several lemmata. 
Recall the definition of admissibily 
of cubes from \S~\ref{sec:absolute geometric presentation theorem} 
and the definition of Koszul cubes associated with graded objects 
from the paragraph before the proof of Theorem~\ref{thm:motive of grF1}.

\begin{para}
\label{lem:admissibility}
{\bf Lemma.}\ \ 
{\it
Let $\cA$ be an abelian category and 
let $n$ and 
$k$ be non-negative integers. 
Then the cube $\kos(x[t_1,\cdots,t_n](-k))$ is admissible. 
}
\end{para}

\begin{proof}
First notice that $\Tot\kos(x[t_1,\cdots,t_n](-k))$ is $0$-spherical 
by \cite[Proposition~4.23 (1)]{MS13}. 
We proceed by induction on $n$. 
For $n=1$, assertion is trivial. 
For $n>1$, notice that a face of $\kos(x[t_1,\cdots,t_n](-k))$ 
is isomorphic to 
$\kos(x[t_1,\cdots,t_{n-1}](-k-j))$ for $j=0$ or $j=1$. 
Therefore by the inductive hypothesis, 
it is admissible. 
Hence by \cite[Corollary~3.15]{Moc13a}, 
$\kos(x[t_1,\cdots,t_n](-k))$ is admissible. 
\end{proof}

\begin{para}
\label{lem:1-s}
{\bf Lemma.}\ \ 
{\it
Let $A$ be the $\bbF_1[t_1,\cdots,t_n]$ the $n$th polynomial 
$\bbN$-grading ring 
over $\bbF_1$. 
Then the functor $T_0(-\otimes_AA[t])\colon\Coh\Nil_{(1)}\gr_A\cA\to 
\Coh\gr_{A[t]}\cA$ is exact and induced map on 
nilpotent invariant dash motives makes the diagram below 
commutative
\begin{equation}
\label{eq:1-s}
\xymatrix{
M'_{\nilp}(\Nil_{(1)}\gr_A\cA) \ar[r]^{M'_{\nilp}(T_0(-\otimes_AA[t]))} \ar[d]_{\wr} & 
M'_{\nilp}(\gr_{A[t]}\cA) \ar[d]^{\wr}\\
\bigoplus_{i\geq 0}M'_{\nilp}(\cA)s^i \ar[r]_{1-s} & 
\bigoplus_{i\geq 0}M'_{\nilp}(\cA)s^i
}
\end{equation}
where the vertical morphisms are isomorphisms in 
Theorem~\ref{thm:motive of grF1}. 
}
\end{para}

\begin{proof}
The functor 
$T_0(-\otimes_AA[t])\colon\Coh\Nil_{(1)}\gr_A\cA \to \Coh\gr_{A[t]}\cA$ 
is exact by Lemma~\ref{lem:admissibility}. 
For an object $x$ in $\Nil_{(1)}\gr_A\cA$, 
there exists an exact sequence
$$x[t](-1)\overset{t}{\rinf} x[t]\rdef T_0(x[t])$$
in $\Coh\gr_{A[t]}\cA$. 
Therefore by additivity theorem and Theorem~\ref{thm:motive of grF1}, 
it turns out that the square $\mathrm{(\ref{eq:1-s})}$ is commutative. 
\end{proof}

In the next lemma, to make the statement simplify, we 
utilize the conventions of affine varieties 
over categories 
in Definition~\ref{df:affine varieties}. 

\begin{para}
\label{lem:split exact}
{\bf Lemma.}\ \ 
{\it
Let $\cB$ be an additive category with 
countable coproduct and $x$ an object 
in $\cB$. 
We consider the object $x[t]$ in $\Spec_{\cB}\bbF_1[t]$ and 
we define $q_x\colon x[t] \to x[t]$, 
$\Diamond_x\colon x[t]\rdef xt^0=x$ and 
$\displaystyle{\pi_x\colon x[t]\rdef \bigoplus_{i=0}^nxt^i}$ 
to be morphisms in $\Spec_{\cB}\bbF_1[t]$
by setting $\displaystyle{q_x:=\begin{pmatrix}
0 & -\id & -\id & -\id & \cdots\\
0 & 0 & -\id & -\id & \cdots\\
0 & 0 & 0 & -\id & \cdots\\
\vdots & \vdots & \vdots & \vdots & \ddots
\end{pmatrix}}$, 
$\Diamond_x=\begin{pmatrix}\id & 0 & 0 & \cdots \end{pmatrix}$ 
and $\displaystyle{\pi_x:=\begin{pmatrix}\Diamond_xq_x^n\\ 
\Diamond_xq_x^{n-1}\\
\vdots\\
\Diamond_x
 \end{pmatrix} }$. 
Then the sequence
\begin{equation}
\label{eq:split exact seq}
x[t]\overset{{(1-t)}^{n+1}}{\rinf}x[t]\overset{\pi_x}{\rdef}\bigoplus_{i=0}^nxt^i
\end{equation}
where $xt^j$ is a just a copy of $x$ 
is a split exact sequence in $\Spec_{\cB}\bbF_1[t]$.
}
\end{para}

\begin{proof}
We set $i_x\colon x\onto{\id_x}x=xt^0\rinc x[t]$ and 
$\displaystyle{j_x\colon \bigoplus_{i=0}^nxt^i\to x[t]}$ by the formula 
$$\displaystyle{j_x=\begin{pmatrix}{(1-t)}^ni_x & {(1-t)}^{n-1}i_x & \cdots 
& i_x \end{pmatrix}}.$$
Then we can check the equalities 
$\pi_x{(1-t)}^{n+1}=0$, $q_x^{n+1}j_x=0$,  
$q_x^{n+1}(1-t)^{n+1}=\id_{x[t]}$, $\pi_xj_x=\id_{\bigoplus_{i=0}^nxt^i}$ and 
${(1-t)}^{n+1}q_x^{n+1}+j_x\pi_x=\id_{x[t]}$. 
This means that the sequence $\mathrm{(\ref{eq:split exact seq})}$ 
is a split exact sequence.
\end{proof}

\begin{proof}[Proof of Corollary~\ref{cor:nilp inv motives of projective spaces}]
We define 
$T\colon \Coh\gr_{\bbF_1}\cA\to\Coh\gr_{\bbF_1[t_0,\cdots,t_n]\cA}$ 
to be a functor by sending an object ${(x_i)}_{i\geq 0}$ to 
$\displaystyle{\bigoplus_{i\geq 0}T_0(x_i[t_0,\cdots,t_n](-i))}$. 
Notice that for an object $ $ in $\Coh\gr_{\bbF_1}\cA$, 
there exists an integer $n>0$ such that $x_i=0$ for $i\geq n$. 
By Lemma~\ref{lem:admissibility}, 
the functor $T$ is exact and since 
$\Coh\Nil_{(1)}\gr_{\bbF_1[t_0,\cdots,t_n]}\cA$ 
is contained in 
$\Coh\gr_{\bbF_1[t_0,\cdots,t_n]}^{V(\ft)}\cA$ 
where $\ft=\{t_i\}_{0\leq i\leq n}$, 
it turns out that the inclusion 
$\Coh\gr_{\bbF_1}\cA\rinc 
\Coh\Nil_{(1)}\gr_{\bbF_1[t_0,\cdots,t_n]}\cA$ 
is a nilpotent immersion 
by Proposition~\ref{prop:fundamental properties of Gr} 
$\mathrm{(3)}$ and $\mathrm{(4)}$. 
Thus it induces an isomorphism of nilpotent invariant dash motives 
$M'_{\nilp}(\Coh\gr_{\bbF_1}\cA)\isoto M'_{\nilp}(\Nil_{(1)}\gr_{\bbF_1[t_0,\cdots,t_n]}\cA)$. 
We will show that this map makes the square \textbf{I} below 
commutative
\begin{equation}
\label{eq:key ladder}
{\scriptscriptstyle{
\xymatrix{
\bigoplus_{i\geq 0}M'_{\nilp}(\cA)s^i \ar[r]^{{(1-s)}^{i+1}} \ar[d]_{\wr}
 \ar@{}[dr]|{\textbf{I}} & 
\bigoplus_{i\geq 0}M'_{\nilp}(\cA)s^i \ar[r] \ar[d]_{\wr} 
\ar@{}[dr]|{\textbf{II}} & 
\bigoplus_{i=0}^nM_{\nilp}'(\cA)a^i \ar[r]^{\!\!\!\!0} \ar[d] 
\ar@{}[dr]|{\textbf{III}} & 
\bigoplus_{i\geq 0}\Sigma M'_{\nilp}(\cA)s^i \ar[d]^{\wr}\\
M'_{\nilp}(\Nil_{(1)}\gr_{A}\cA) \ar[r] &
M'_{\nilp}(\gr_{A}\cA) \ar[r] & 
M'_{\nilp}(\bbP^n_{\cA}) \ar[r]_{\!\!\!\!\!\!\!\!\!\!\!\!\!\!\!\!\!\!\!\!\!\!\partial} & 
\Sigma M'_{\nilp}(\Nil_{(1)}\gr_{A}\cA)
}}}
\end{equation}
where for simplicity we set $A=\bbF_1[t_0,\cdots,t_n]$. 
We consider the following commutative diagram
$$\xymatrix{
\Coh\gr_{\bbF_1}\cA \ar[rd]^{\ \ \ T_0(-\otimes_{\bbF_1}\bbF_1[t_0,\cdots,t_n])} 
\ar[d]_{T_0(-\otimes_{\bbF_1}\bbF_1[t_0,\cdots,t_n])}\\
\Coh\Nil_{(1)}\gr_{\bbF_1[t_0,\cdots,t_{n-1}]}\cA \ar[r] 
\ar[d]_{T_0(-\otimes_{\bbF_1[t_0,\cdots,t_{n-1}]}\bbF_1[t_0,\cdots,t_n])} & 
\Coh\gr_{\bbF_1[t_0,\cdots,t_{n-1}]}\cA 
\ar[d]^{T_0(-\otimes_{\bbF_1[t_0,\cdots,t_{n-1}]}\bbF_1[t_0,\cdots,t_n])}\\
\Coh\Nil_{(1)}\gr_{\bbF_1[t_0,\cdots,t_{n}]}\cA \ar[r] & 
\Coh\gr_{\bbF_1[t_0,\cdots,t_{n}]}\cA.
}$$
For an object $x$ in $\cA$, 
Inspection shows that an equality 
$$T_0(x[t_0,\cdots,t_n](-k))=
T_0(
T_0(x[t_0,\cdots,t_{n-1}])(-k)\otimes_{\bbF_1[t_0,\cdots,t_{n-1}]}\bbF_1[t_0,\cdots,t_n]).$$
Thus by induction on $n$ and Lemma~\ref{lem:1-s}, 
it turns out that the square \textbf{I} is commutative. 

Next since both the top and the bottom lines are distinguished 
triangles by Lemma~\ref{lem:split exact} 
and Lemma~\ref{lem:Nil-> Gr-> proj} 
respectively, 
there exists a morphism $\displaystyle{\bigoplus_{i=0}^nM'_{\nilp}(\cA)s^i 
\to M'_{\nilp}(\bbP_{\cA}^n)}$ which makes the square \textbf{II} and 
\textbf{III} commutative and by five lemmam this map is an isomorphism. 
Thus it turns out that the bottom line is also a split exact sequence. 
Namley the map $\partial$ in the diagram $\mathrm{(\ref{eq:key ladder})}$ is trivial. 
Thus we learn that the map 
$\displaystyle{\bigoplus_{i=0}^nM'_{\nilp}(\cA)s^i\to M'_{\nilp}(\bbP^n_{\cA}) }$ 
induced from the functor 
$\displaystyle{\bigoplus_{i=0}^n\Coh\cA\to \bbP^n_{\cA}}$ 
which sends an object 
$(x_i)_{0\leq i\leq n}$ to 
$\displaystyle{\bigoplus_{i=0}^nx_i[t_0,\cdots,t_n](-i)}$ also makes 
the square \textbf{II} and \textbf{III} commutative and by 
five lemma again, this map is also an isomorphism. 
We complete the proof.  
\end{proof}

Similarly we define the notion of affine varieties over categories. 

\begin{para}
\label{df:affine varieties}
{\bf Definition (Affine varieties).}\ \ 
Let $A$ be a commutative $B$-algebra and 
let $\cA$ be a noetherian abelian $B$-category. 
We regard $A$ as a $B$-category by setting 
$\Ob A=\{\ast\}$ and $\Hom_A(\ast,\ast):=A$. 
Then we define 
$\Spec_{\cA} A $ to be a $B$-category 
by setting $\Spec_{\cA} A:=\HOM_B(A,\cA)$ 
the category of $B$-functors from $A$ to $\cA$ and natural transformations 
and call it 
an {\it affine scheme associated with $A$ 
over $\cA$}. 
Namely an object in $\Spec_{\cA} A$ can be regarded as a pair of 
an object $x$ in $\cA$ and a family of endomorphisms 
$\{x_f\colon x\to x\}_{f\in A}$ indexed 
by elements in $A$ such that endomorphisms satisfy equalities arose 
from relations among elements in $A$. 
We sometimes abbreviate the morphism $x_f\colon x\to x$ to $f$. 

For a family of elements $\ff=\{f_k\}_{1\leq k\leq r}$ of $A$, 
let $\Spec_{\cA}^{V(\ff)}A $ be a full subcategory of $\Spec_{\cA} A$ consisting of 
those objects $x$ such that there exists a positive integer $N>0$ such that 
$x_{f_k}^N=0$ for all $1\leq k\leq r$ and we 
denote the full subcategory of $\Spec_{\cA}^{V(\ff)}A$ 
consisting of those objects $x$ such that 
$x_{f_k}=0$ for all $1\leq k\leq r$ by $\Spec_{\cA,\red}^{V(\ff)}A$. 

By convention, we write $\Spec_{\cA}\bbF_1$ for $\cA$ and for any positive 
integer $n$, 
we set $\bbA^n_{\cA}:=\Spec_{\cA}B[t_1,\cdots,t_n]$ 
and call it the {\it $n$th affine space over $\cA$}. 

By virtue of 
Proposition~\ref{prop:fundamental properties of affine varieties} 
$\mathrm{(9)}$ below, 
we can regard $\Spec A$ as a functor from 
the (large) category of locally noetherian abelian category and cocontinuous 
functors to itself.

Let $C$ be a $\bbN$-grading commutative $B$-algebra. 
Then there exists a forgetting grading functor 
$F_C\colon \gr_C\cA \to \Spec_{\cA}C$ which sends an object $x$ 
in $\gr_C \cA$ to a pair $\displaystyle{\left(\bigoplus_{n\geq 0}x_n, 
\{F_A(x)_f \}_{f\in C} \right)}$. 
Here for a homogeneous element $f\in C$ of degree $s$, 
we set $\displaystyle{{F_A(x)}_f:=\bigoplus_{k\geq 0}f
\colon x_k\to x_{k+s}\colon 
\bigoplus_{n\geq 0}x_n \to \bigoplus_{n\geq 0}x_n }$. 
For a general element $f\in C$, we denote it by the summation of 
homogeneous elements $\displaystyle{f=\sum_{k=1}^r f_k }$, $f_k\in A_k$ and 
we set $\displaystyle{{F_A(x)}_f:=\sum_{k=1}^r {F_A(x)}_{f_k} }$. 
Since $\gr_A\cA$ is locally noetherian, 
the operation $\bigoplus$ is exact and 
therefore $F_A$ is a faithful and exact functor.
\end{para}

\begin{para}
\label{ex:affine spectrum}
{\bf Examples (Affine varieties).}\ \ 
Let $\cA$ be a locally noetherian abelian $B$-category. 
\begin{enumerate}
\enumidef
\item
Let $X$ be a noetherian scheme and $n$ be a positive integer. 
Then we have the canonical equivalence of categories 
$\bbA^n_{\Qcoh_X}\isoto \Qcoh_{\bbA^n_X}$.

\item
Let $A$ and $C$ be a commutative $B$-algebras. 
Then we have the canonical equivalence of categories 
$\Spec_{\Mod_A}C\isoto \Mod_{A\otimes_B C}$. 
\end{enumerate}
\end{para}

\begin{para}
\label{prop:fundamental properties of affine varieties}
{\bf Proposition (Fundamental properties of affine varieties).}\ \ 
{\it Let $\cA$ be a category and $A$ a commutative $B$-algebra. Then
\begin{enumerate}
\enumidef
\item
We can calculate a {\rm{(}}co{\rm{)}}limit 
in $\Spec_{\cA}A$ by 
term-wise {\rm{(}}co{\rm{)}}limit in $\cA$. 
In particular, if 
$\cA$ is additive {\rm(}resp. abelian, 
closed under small filtered colimits{\rm)} 
then $\Spec_{\cA}A$ is also additive {\rm(}resp. abelian, closed under 
small filtered colimits{\rm)}.

\mn
Moreover assume that $\cA$ is abelian and let $\ff=\{f_k\}_{1\leq k\leq r}$ 
be a family of homogeneous elements of $A$. 
We write $\ff A$ for the ideal of $A$ spanned by $\ff$. Then

\item
$\Spec_{\cA}^{V(\ff)}A$ is a Serre subcategory of $\Spec_{\cA}A$. 

\item
$\Spec_{\cA,\red}^{V(\ff)}A$ is a topologizing subcategory 
of $\Spec_{\cA}^{V(\ff)}A$ and the inclusion functor 
$\Spec_{\cA,\red}^{V(\ff)}A\rinc \Spec_{\cA}^{V(\ff)}A$ is 
a nilpotent immersion. 

\item
The canonical map $A\rdef A/\ff A$ induces an equivalence of categories 
$\Spec_{\cA} A/\ff A\isoto \Spec_{\cA,\red}^{V(\ff)}A$. 

\item
The fully faithful embedding $\Spec_{\cA}A/\ff A\rinc \Spec_{\cA}A$ indcued from the canonical map 
$A\rdef A/\ff A$ admits a left adjoint functor 
$-\otimes_A A/\ff A\colon \Spec_{\cA}A \to \Spec_{\cA}A/\ff A$.

\mn
Moreover assume that $A=B[t_1,\cdots,t_n]/(g_1,\cdots,g_t)$ 
where $B[t_1,\cdots,t_n]$ is the $n$th polynomial ring over $B$ and 
$g_i$ is a polynomial in $B[t_1,\cdots,t_n]$ 
for $1\leq i\leq t$. Then

\item
The forgetful functor $U_A\colon \Spec_{\cA}A\to \cA$ which 
sends an object $(x,\{x_f\}_{f\in A})$ in $\Spec_{\cA}A$ to $x$ admits 
a left adjoint functor $-\otimes_B A\colon \cA \to \Spec_{\cA}A$. 

\item
If $x$ in $\cA$ is a noetherian object, 
then $x\otimes_B A$ is also 
a noetherian object in $\Spec_{\cA}A$. 

\item
For any object $x$ in $\cA$, 
there exists a canonical epimorphism 
$U_A(x)\otimes_B A\rdef x$.

\item
If $\cA$ admits a family of generators 
$\{u_{\lambda}\}_{\lambda\in\Lambda}$ indexed by non-empty set $\Lambda$, 
then 
$\Spec_{\cA} A$ has a family of generators 
$\{u_{\lambda}\otimes_B A\}_{\lambda\in\Lambda}$. 
In particular if $\cA$ is Grothendieck 
{\rm (}resp. locally noetherian{\rm )}, 
then $\Spec_{\cA}A$ is also.
\end{enumerate}
}
\end{para}

\begin{proof}
$\mathrm{(1)}$, $\mathrm{(2)}$, $\mathrm{(3)}$ and $\mathrm{(4)}$ 
are straightforward. 

\sn
$\mathrm{(5)}$ 
As in the proof of 
Proposition~\ref{prop:fundamental properties of Gr} $\mathrm{(5)}$, 
for an object $x$ in $\Spec_{\cA} A$, 
we define $\ff x$ to be a subobject by setting 
$\displaystyle{\ff x:=\sum_{k=1}^r \im (f_k\colon x\to x)}$ 
and we set $x\otimes_A A/\ff A:=x/\ff x$. 
Then we can show that the association $x\mapsto x\otimes_AA/\ff A$ 
gives a left adjoint functor of the fully faithful embedding 
$\Spec_{\cA}A/\ff A\rinc \Spec_{\cA}A$.

\sn
$\mathrm{(6)}$ 
Since $U_A\colon\Spec_{\cA}A \to \cA$ factors through 
$$\Spec_{\cA}A \rinc \Spec_{\cA}B[t_1,\cdots,t_n]\to \cA,$$
by $\mathrm{(5)}$, we shall only prove the case where 
$A=B[t_1,\cdots,t_n]$. 
In this case the functor 
$-\otimes_B B[t_1,\cdots,t_n]\colon \cA\to \Spec_{\cA} A$ 
is just a composition of 
$\cA \onto{-\otimes_BB[t_1,\cdots,t_n]} \gr_{B[t_1,\cdots,t_n]}\cA \onto{F_{B[t_1,\cdots,t_n]}} \Spec_{\cA} A$.

\sn
$\mathrm{(7)}$ 
As in the proof of 
Proposition~\ref{prop:fundamental properties of Gr} $\mathrm{(8)}$, 
we can reduce to the case where $A=B[t_1,\cdots,t_n]$ and in this case, 
assertion is proven in \cite[9.10 b]{Sch00}.

\sn
$\mathrm{(8)}$ 
Let $x$ be an object in $\Spec_{\cA}A$. 
Then there exists the canonical projection 
$$\displaystyle{U_A(x)[t_1,\cdots,t_n]=\bigoplus_{(i_1,\cdots,i_n)\in\bbN} xt_1^{i_1}\cdots t_n^{i_n}
\rdef xt_1^0\cdots t_n^0=x}$$
and this morphism factors through 
$U_A(x)[t_1,\cdots,t_n]\rdef U_A(x)\otimes_B A \rdef x$. 
The last morphism is a desired epimorphism.

\sn
$\mathrm{(9)}$ 
For a non-zero morphism $x\onto{a} y$ in $\Spec_{\cA} A$, 
since $F_A(a)$ is a non-zero morphism, there exists 
a morphism $b\colon u_\lambda \to F_A(x)$ such that $ab\neq 0$. 
Then the compositions 
$u_{\lambda}\otimes_B A\onto{b\otimes_B A} F_A(x)\otimes_B A \rdef x
\onto{a} y$ is a non-zero morphism. 

\end{proof}

For an $\bbN$-grading commutative $B$-algebra $A$ and 
a locally noetherian abelian $B$-category, 
recall the functor $F_A\colon \gr_A\cA\to \Spec_{\cA}A$ 
from Definition~\ref{df:affine varieties}. 
The functor sends an object $x$ in $\gr_A\cA$ to 
a pair 
$\displaystyle{\left (\bigoplus_{n\geq 0} x_n, \left \{ 
{F_A(x)}_f
\right \}_{f\in A} \right ) }$ in $\Spec_{\cA}A$.

\begin{para}
\label{prop:about FA}
{\bf Poposition.}\ \ 
{\it
Let $A$ be an $\bbN$-grading commutative $B$-algebra. Then
\begin{enumerate}
\enumidef
\item
For an object $x$ in $\gr_A\cA$ and non-negative integer $k\geq 0$, 
$F_A(x(-k))=F_A(x)$. 

\item
For a family of homogeneous elements $\ff=\{f_k\}_{1\leq k\leq r}$ in 
$A$, we have the equality 
$F_A(\ff x)=F_A(\ff)F_A(x)$ where $F_A(\ff):=\{F_A(f_k)\}_{1\leq k\leq r}$. 

\item 
Let $y$ and $z$ be a pair of subobjects of $x$, then we 
have the equiality $F_A(y\cap z)=F_A(y)\cap F_A(z)$.

\mn
Moreover assume that there exists a finite family of homogeneous elements 
$\{g_i \}_{1\leq i\leq r}$ of $B[t_1,\cdots,t_n]$ such that 
$A=B[t_1,\cdots,t_n]/(g_1,\cdots,g_r)$. 

\item
For an object $x$ in $\cA$, we have the equality 
$F_A(x\otimes_B A)=x\otimes_BA$. 

\item
$F_A$ induces the functor $\Coh\gr_A\cA\to\Coh\Spec_{\cA}A$. 
Namely $F_A$ sends a coherent object $x$ in $\gr_A\cA$ to a 
coherent object $F_A(x)$ in $\Spec_{\cA}A$. 

\end{enumerate}
}
\end{para}

\begin{proof}
$\mathrm{(1)}$ 
$\displaystyle{F_A(x(-k))=\bigoplus_{n\geq k}x_{n-k}=\bigoplus_{n\geq 0}x_n=F_A(x) }$.

\sn
$\mathrm{(2)}$ 
\begin{eqnarray*}
F_A(\ff x)=\sum_{k=1}^r\im\left (\bigoplus_{n\geq 0}f_k 
\colon \bigoplus_{n\geq 0}x_n \to \bigoplus_{n\ge 0}x_n \right )\\
=\sum_{k=1}^r\bigoplus_{n\geq 0}\im (f_k\colon x_{n-\deg f_k}\to x_n )\\
=\bigoplus_{n\geq 0}\sum_{k=1}^r \im (f_k\colon x_{n-\deg f_{k}}\to x_n)
=F_A(\ff)F_A(x).
\end{eqnarray*}

\sn
$\mathrm{(3)}$ 
$\displaystyle{F_A(y\cap z)=\bigoplus_{n\geq 0}(y_n\cap z_n)
=\left (\bigoplus_{n\geq 0}y_n  \right ) \cap \left (\bigoplus_{n\geq 0}z_n \right )=F_A(y)\cap F_A(z)  }$.

\sn
$\mathrm{(4)}$ By exactness of $F_A$ and $\mathrm{(2)}$, 
we have the equalities 
\begin{eqnarray*}
F_A(x\otimes_BA)=
F_A(x\otimes_BB[t_1,\cdots,t_n]\otimes_{B[t_1,\cdots,t_n]}A)\\
=F_{B[t_1,\cdots,t_n]}(x\otimes_BB[t_1,\cdots,t_n])
\otimes_{B[t_1,\cdots,t_n]} A\\
=x\otimes_BB[t_1,\cdots,t_n]
\otimes_{B[t_1,\cdots,t_n]}A=x\otimes_BA.
\end{eqnarray*}

\sn
$\mathrm{(5)}$ 
Assume that $x$ is noetherian, then there exists an integer $k>0$ 
and an epimorphism 
$\displaystyle{\bigoplus_{k=0}^mx_k\otimes_BA(-k) \rdef x}$. 
Then by exactness of $F_A$ and $\mathrm{(1)}$, $\mathrm{(4)}$, 
we obtain the epimorphism  
$\displaystyle{\bigoplus_{k=0}^mx_k\otimes_BA \rdef F_A(x)}$. 
Thus since $F_A(x)$ is a quotient of finite direct sum of 
noetherian objects, $F_A(x)$ is also a noetherian object.
\end{proof}

\begin{para}
\label{prop:fundamental properties of dash motives}
{\bf Proposition.}\ \ 
{\it
Let $\cA$ be a locally noetherian abelian $B$-category and 
let $A$ be a $\bbN$-grading commutative $B$-algebra 
$\ff=\{f_i\}_{1\leq i\leq s}$ a 
finite non-empty family of homogeneous elements in $A$ 
and 
let $C$ be a commutative $B$-algebra and 
$\fg=\{g_i \}_{1\leq i\leq t}$ a 
finite non-empty family of elements in $C$. 
Then
\begin{enumerate}
\enumidef
\item 
{\bf (Localization distinguished triangles).}
For $\#\in\{\loc,\nilp \}$, there exists distinguished triangles 
of dash motives 
$${\scriptstyle{{M'}_{\nilp}^{V(\fg)}(\Spec_{\cA} C)\to M'_{\nilp}(\Spec_{\cA} C) 
\to M'_{\nilp}\left(\Spec_{\cA}C/\Spec^{V(\fg)}_{\cA}C \right) \to \Sigma {M'}_{\nilp}^{V(\fg)}(\Spec_{\cA} C),}}$$
$${\scriptstyle{{M'}_{\nilp}^{V(\ff)}(\proj_{\cA} A)\to M'_{\nilp}(\proj_{\cA} A) 
\to M'_{\nilp}\left(\proj_{\cA}A/\proj^{V(\ff)}_{\cA}A \right) \to 
\Sigma {M'}_{\nilp}^{V(\ff)}(\proj_{\cA} A)}} $$
in $\Mot_{\dg}^{\#}$. 

\item
{\bf (Purity).}\ \ 
The canonical inclusions 
$\Spec_{\cA}C/\fg C\rinc \Spec_{\cA}^{V(\fg)}C$ and 
$\proj_{\cA}A/\ff A\rinc \proj_{\cA}^{V(\ff)}A$ induce 
isomorphisms of nilpotent invariant dash motives 
$M'_{\nilp}(\Spec_{\cA} C/\fg C)
\isoto{M'}_{\nilp}^{V(\fg)}(\Spec_{\cA}C)$ and 
$M'_{\nilp}(\proj_{\cA} A/\ff A)
\isoto{M'}_{\nilp}^{V(\ff)}(\proj_{\cA}A) $.

\end{enumerate}
}
\end{para}

\begin{proof}
$\mathrm{(1)}$ 
We will only give a proof for the projective case. 
A proof of the affine case is similar. 
For simplicity, we set $\cB=\Coh\gr_A\cA$ and $\cC=\Coh\gr_A^{V(\ff)}\cA$. 
Then there exists a $3\times 3$ commutative diagrams 
$$
{\footnotesize{\xymatrix{
\Nil_{(1)}\cB \ar[r] \ar[d] & \Nil_{(1)}\cC \ar[r] \ar[d] & 
\Nil_{(1)}\cC/\Nil_{(1)}\cB \ar[d]\\
\cB \ar[r] \ar[d] & \cC \ar[r] \ar[d] & \cC/\cB \ar[d]\\
\Coh\proj_{\cA}^{V(\ff)} A \ar[r] & 
\Coh\proj_{\cA} A \ar[r] &
\Coh \left (\proj_{\cA}A/\proj_{\cA}^{V(\ff)}A \right ).
}}}
$$
What we need to prove is that 
$\cB$ and $\Nil_{(1)}\cB$ are special filtering 
in $\cC$ and $\Nil_{(1)}\cC$ 
respectively. 
Then assertion follows from $3\times 3$-lemma for distinguished 
triangles. (\cf \cite[Exercise 10.2.6.]{Wei94}.)

Let $x\rinf y$ be a monomorphism in $\cC$ with $x$ in $\cB$. 
Thus there exists a non-negative integer $m\geq 0$ such that $\ff^mx=0$. 
Then by Artin-Rees lemma~\ref{lem:char of stability} below, 
we can show that there exists a non-negative integer $n_0\geq 0$ such that 
$\ff^ny\cap x=\ff^{n-n_0}(\ff^{n_0}y\cap x)=0$ for $n\geq n_0+m$. 
Hence for $n\geq n_0+m$, 
the composition $x \rinf y \rdef y/\ff^ny$ is a monomorphim 
and $y/\ff^ny$ is in $\cB$. 
Thus $\cB$ is right special filtering in $\cC$. 
A proof of that $\Nil_{(1)}\cB$ is 
right special filtering in $\Nil_{(1)}\cC$ is similar. 

\sn
$\mathrm{(2)}$ 
For affine case, 
it is a direct consequence of 
Proposition~\ref{prop:fundamental properties of affine varieties} 
$\mathrm{(3)}$ and 
nilpotent invariance of $M_{\nilp}$. 

For projective case, we consider the commutative diagram of 
distinguished triangles 
$${\footnotesize{
\xymatrix{
M'_{\nilp}(\Nil_{(1)}\gr_{A/\ff A}\cA)\ar[r] \ar[d]^{\textbf{I}} &
M'_{\nilp}(\gr_{A/\ff A}\cA) \ar[r] \ar[d]_{\textbf{II}} & 
M'_{\nilp}(\proj_{A/\ff A}\cA) \ar[r]^{\ \ \ \ \ \ \ \ \ +1} \ar[d]_{\textbf{III}}& \\ 
M'_{\nilp}(\Nil_{(1)}\gr_{A}^{V(\ff)}\cA)\ar[r] &
M'_{\nilp}(\gr_{A}^{V(\ff)}\cA) \ar[r] & 
{M'}_{\nilp}^{V(\ff)}(\proj_{A}\cA) \ar[r]_{\ \ \ \ \ \ \ \ \ +1} &.
}}}
$$
As in Proposition~\ref{prop:fundamental properties of Gr} $\mathrm{(3)}$, 
we can also show that 
$\Nil_{(1)}\gr_{A/\ff A}\cA$ is topologizing subcategory of 
$\Nil_{(1)}\gr_{A}^{V(\ff)}\cA $ and 
the inclusion functor 
$\Nil_{(1)}\gr_{A/\ff A}\cA\rinc\Nil_{(1)}\gr_{A}^{V(\ff)}\cA$ 
is a nilpotent immersion. 
Thus by nilpotent invariance of $M_{\nilp}$ the morphisms $\textbf{I}$, 
$\Sigma\textbf{I}$ 
and $\textbf{II}$ are isomorphisms and by five lemma, 
the morphism $\textbf{III}$ is also isomorphism. 

\end{proof}

To state an abstract Artin-Rees lemma~\ref{lem:char of stability} below, 
we will introduce some notations 
from \cite[\S 2.5]{MS13} with slight generaization. 
Let $A$ be a commutative $B$-algebra of finite type over $B$ and let 
$\cA$ be a locally noetherian abelian $B$-category and 
$\ff=\{f_i\}_{1\leq i\leq r}$ be 
a non-empty finite family of elements in $A$. 
For an object $x$ in $\Spec_{\cA}A$, recall that 
we denote the subobject $\displaystyle{\sum_{i=1}^r\im(f_i\colon x\to x)}$ 
of $x$ by $\ff x$ and inductively we set ${\ff}^{n+1}x:=\ff({\ff}^nx)$ for 
a non-negative integer $n\geq 0$. 

For a decreasing filtration 
$\fx=\{x_n\}_{n\geq 0}$ of $x$, 
$x=x_0\linf x_1\linf x_2\linf \cdots \linf x_n \linf \cdots$ 
is a {\it $\ff$-filtration} if 
$\ff x_n\subset x_{n+1}$ for any $n\geq 0$. 
A $\ff$-filtration $\fx=\{x_n\}_{n\geq 0}$ is {\it stable} 
if there exists an integer $n_0\geq 0$ such that 
$\ff x_n= x_{n+1}$ for any $n\geq n_0$.

For a $\ff$-filtration $\fx=\{x_n\}_{n\geq 0}$ of $x$, 
We define $\Bl_{\ff}\fx$ 
to be an object in $\Spec_{\cA}A$ by setting 
$$\displaystyle{\Bl_{\ff}\fx:=\left (\bigoplus_{n\geq 0} x_n, 
\left \{\bigoplus_{n\geq 0}{x_n}_f \right\}_{f\in A} \right )}.$$
We call $\Bl_{\ff}\fx$ a {\it blowing up object of $\fx$ along $\ff$}. 
For each non-negative integer $n$ and each  $1\leq k\leq r$, 
the morphisms $f^p_k\colon x_n \to x_{n+p}$ 
for $p>1$ induce a morphism 
$\eta_{f_k}^n\colon x_n\to \Bl_{\ff}\fx$ 
in $\Spec_{\cA} A$.

\begin{para}
\label{lem:char of stability}
{\bf Lemma.}\ \ 
{\it 
Let $A$ be a commutative $B$-algebra of finite type over $B$ and 
$\ff=\{f_i\}_{1\leq i\leq r}$ be 
a non-empty finite family of elements in $A$ and 
let $C$ be an $\bbN$-grading commutative $B$-algebra of finite type over $B$ and $\fg=\{g_j\}_{1\leq j\leq s}$ be 
a non-empty finite family of homogeneous elements in $C$ and 
let $\cA$ be a locally noetherian abelian category. 
Then
\begin{enumerate}
\enumidef
\item
{\bf (Characterization of stability).}\ \ 
For an object $x$ in $\Coh\Spec_{\cA}A$ and 
$\ff$-filtration $\fx:=\{x_n\}_{n\geq 0}$ of $x$, 
the following conditions are equivalent:
\begin{enumerate}
\enumiidef 
\item
$\fx$ is stable.

\item 
There exists an integer $m\geq 0$ such that the canonical morphism 
induced by $\eta_{f_i}^k$ {\rm (}$0\leq k\leq m$, $1\leq i\leq r${\rm )}, 
$\displaystyle{\bigoplus_{i=1}^r\bigoplus_{k=0}^mx_k \to \Bl_{\ff}\fx}$ 
is an epimorphism.

\item
$\Bl_{\ff}\fx$ is an object in $\Coh\Spec_{\cA}A$, 
namely a noetherian object in $\Spec_{\cA} A$.
\end{enumerate}

\item 
{\bf (Abstract Artin-Rees lemma).}
\begin{enumerate}
\enumiidef
\item
Let $x$ be an object in $\Coh\Spec_{\cA} A$ and $y$ be 
a subobject of $x$. 
Then there exists a non-negative integer $n_0\geq 0$ such that 
$\ff^n x\cap y=\ff^{n-n_0}\left(\ff^{n_0}x\cap y \right )$ 
for any $n\geq n_0$.

\item
Let $x$ be an object in $\Coh\gr_C \cA$ and $y$ be a subobject of $x$. 
Then there exists a non-negative integer $n_0\geq 0$ 
such that 
$\fg^n x\cap y=\fg^{n-n_0}\left(\fg^{n_0}x\cap y \right )$ 
for any $n\geq n_0$.
\end{enumerate}
\end{enumerate}
}
\end{para}

\begin{proof}
$\mathrm{(1)}$ 
We assume that there exists an integer $m\geq 0$ such that 
$\ff x_n=x_{n+1}$ for any $n\geq m$. 
Then obviously the canonical morphism 
$\displaystyle{\bigoplus_{i=1}^r\bigoplus_{k=0}^mx_k \to \Bl_{\ff}\fx}$ 
is an epimorphism. 

\sn
Next assume the condition $\mathrm{(ii)}$. 
Since $\Bl_{\ff}\fx$ is a quotient of finite direct sum of 
noetherian objects 
in $\Spec_{\cA}A$, 
$\Bl_{\ff}\fx$ is a noetherian object. 

\sn
Finally we assume that $\Bl_{\ff}\fx$ is noetherian. 
We set $\displaystyle{z_m\colon =
\im\left (\bigoplus_{i=1}^r\bigoplus_{k=0}^mx_k \to \Bl_{\ff}\fx \right )}$. 
Then the sequence $z_0 \rinf z_1\rinf z_2\rinf \cdots $ is stational. 
Say $z_{n_0}=z_{n_0+1}=\cdots$. 
Then for any $n\geq n_0$, we have 
$$x_{n+1}\subset z_{n+1}\cap x_{n+1}=z_{n_0}\cap x_{n+1} 
\subset \sum_{k=1}^r\sum_{i=0}^{n_0} \im (f_k^{n+1-i}\colon x_i \to x_{n+1})\subset 
\ff x_{n}.$$
Hence $\fx$ is stable.

\sn
$\mathrm{(2)}$ 
$\mathrm{(i)}$ 
Consider the $\ff$-stable filtration 
$\fx=\{{\ff}^nx \}_{n\geq 0}$ of $x$ and 
the induced $\ff$-filtration 
$\fy=\{({\ff}^nx)\cap y \}_{n\geq 0}$ of $y$. 
Then $\Bl_{\ff}\fy$ is a subobject of $\Bl_{\ff}\fx$. 
Since $\Bl_{\ff}\fx$ is noetherian by $\mathrm{(1)}$, 
$\Bl_{\ff}\fy$ is also noetherian and 
by $\mathrm{(1)}$ again, 
we learn that $\fy$ is stable. 
Hence we obtain the result.

\sn
$\mathrm{(ii)}$ 
Let $F_C\colon\Coh\gr_C\cA \to \Coh\Spec_{\cA}C$ 
be the forgetting grading functor. 
(Well-definedness follows 
from Proposition~\ref{prop:about FA} $\mathrm{(5)}$.) 
By applying $F_C(\fg)$, $F_C(x)$ and $F_C(y)$ to 
$\mathrm{(2)}$ $\mathrm{(i)}$ and 
Proposition~\ref{prop:about FA} $\mathrm{(2)}$ and $\mathrm{(3)}$, 
we obtain the equality $F_C(\fg^n x\cap y)=F_C(\fg^{n-n_0}\left(\fg^{n_0}x\cap y \right ))$. 
Since $F_C$ is faithful and exact, we obtain the desired equality. 

\end{proof}

\begin{para}
\label{cor:A^1-homotopy invariance of dash motives}
{\bf Corollary ($\bbA^1$-homotopy invariance).}\ \ 
{\it
Let $\cA$ be a locally noetherian category. 
Then the base change functor 
$-\otimes_{\bbF_1}\bbF_1[t]\colon \cA\to \bbA^1_{\cA}$ 
is exact and 
induces an isomorphism 
$
M'_{\nilp}(\cA)\isoto M'_{\nilp}\left (\bbA^1_{\cA} \right )
$
of nilpotent invariant dash motives. 
}
\end{para}

\begin{proof}
Since $\cA$ is Gothendeick, the operation $\bigoplus$ is exact. 
Hence the functor 
$-\otimes_{\bbF_1}\bbF_1[t]\colon\cA \to \bbA^1_{\cA}$ is exact. 
The functor $\gr_{\bbF_1[t_0,t_1]}\cA\ \to \Spec_{\cA}\bbF_1[t_0]$ which sends an object $x$ to 
$\displaystyle{\left (\underset{t_1}{\colim}\ x_n,\underset{n}{\colim}\ t_0 \right )}$ 
induces an equivalence of categories 
$\Coh\left (\bbP^1_{\cA}/\proj_{\cA}^{V(t_1)}\bbF_1[t_0,t_1] \right )\isoto \Coh\bbA^1_{\cA}$ by \cite[Theorem 5.6]{MS13}. 
Thus replacing $X$ with $\cA$ and setting $n=1$, 
we obtain the commutative diagram 
$\mathrm{(\ref{eq:key commutative diagram})}$ where 
both the top and the bottom lines are distinguished triangles by 
Corollary~\ref{cor:nilp inv motives of projective spaces} 
and 
Proposition~\ref{prop:fundamental properties of dash motives} $\mathrm{(1)}$ 
respectively 
and the morphisms $\textbf{I}$ and $\Sigma\textbf{I}$ 
are isomorphisms by 
Proposition~\ref{prop:fundamental properties of dash motives} $\mathrm{(2)}$. 
Thus we obtain the result by five lemma. 
\end{proof}

\begin{para}
\label{cor:A1-homotopy invariance}
{\bf Corollary (Homotopy invariance and projective bundle formula for dash motives).}\ \ 
{\it Let $X$ be a noetherian scheme. Then
\begin{enumerate}
\enumidef
\item
Let $f\colon P\to X$ be a flat morphism of finite type whose fibers are 
affine spaces. 
Then $f$ induces an isomorphism of 
nilpotent invariant 
dash motives $M'_{\nilp}(P)\isoto M'_{\nilp}(X)$. 
In particular if $X$ and $P$ are regular separated noetherian, then $f$ 
induces an isomorphism of nilpotent invariant dash motives 
$M_{\nilp}(P)\isoto M_{\nilp}(X)$.

\item
Let $\cE$ be a vector bundle of rank $r$ on $X$ and let 
$\bbP(\cE)$ be associated projective bundle. 
Then we have a canonical isomorphism of nilpotent invariant dash 
motives 
$\displaystyle{M'_{\nilp}(\bbP(\cE))\isoto\bigoplus_{k=0}^{r}M'_{\nilp}(X)}$.
\end{enumerate}
}
\end{para}

\begin{proof}
$\mathrm{(1)}$ 
Let $Z$ be a closed subset of $X$, there exists a commutative 
diagram of distinguished triangles of nilpotent invariant motives
$$
\xymatrix{
M'_{\nilp}(Z) \ar[r] \ar[d] & M'_{\nilp}(X) \ar[r] \ar[d] & M'_{\nilp}(X\ssm Z) 
\ar[d] \ar[r] & \Sigma M'_{\nilp}(Z) \ar[d]\\
M'_{\nilp}(P_Z) \ar[r]  & M'_{\nilp}(P) \ar[r] & 
M'_{\nilp}(P_{X\ssm Z}) \ar[r] & \Sigma M'_{\nilp}(P_Z).
}
$$
Then as in the proof of \cite[\S 7 Proposition 4.1.]{Qui73}, 
by utilizing noetherian induction, we can reduce the case where 
$X$ is the spectrum $\Spec k$ of a field $k$.
Then assertion follows form Corollary~\ref{cor:A1-homotopy invariance}.

\sn
$\mathrm{(2)}$ 
Similary we can reduce to the case where $X$ is the spectrum $\Spec k$ of a field $k$. 
Then assertion follows from Corollary~\ref{cor:nilp inv motives of projective spaces}. 

\end{proof}

\medskip
\noindent
SATOSHI MOCHIZUKI\\
{\it{DEPARTMENT OF MATHEMATICS,
CHUO UNIVERSITY,
BUNKYO-KU, TOKYO, JAPAN.}}\\
e-mail:{\tt{mochi@gug.math.chuo-u.ac.jp}}

\end{document}